\documentclass[pdftex]{siamltex}
\usepackage[pdftex]{graphicx}
\usepackage{bmpsize}
\usepackage{bm}
\usepackage{amssymb,amsfonts,amsmath}
\usepackage{color}
\usepackage{boxedminipage}
\usepackage{slashbox}
\usepackage{cancel}
\usepackage[normalem]{ulem}

\newcommand{\refeqn}[1]{(\ref{#1})}
\newcommand{\reffig}[1]{FIG. \ref{#1}}

\title{Topological computation analysis of meteorological time-series data}

\author{
Hidetoshi Morita%
\thanks{
{\em Corresponding author}\,; 
Department of Mathematics,
Kyoto University,
Kyoto 606-8502, Japan
({\tt hmorita@math.kyoto-u.ac.jp})}
\and
Masaru Inatsu%
\thanks{
Faculty of Science, Hokkaido University,
Sapporo, Japan
({\tt inaz@sci.hokudai.ac.jp})}
\and
Hiroshi Kokubu%
\thanks{\,Department of Mathematics,
Kyoto University,
Kyoto 606-8502, Japan
({\tt kokubu@math.kyoto-u.ac.jp})}
}

\begin{document}

\maketitle

\begin{abstract}
A topological computation method (called the MGSTD method) is applied to noisy time-series 
data obtained from meteorological measurement.
This method is based on the idea of the
Morse decomposition, which is a decomposition of the dynamics into invariant sets, called 
the Morse sets, and their gradient-like connections.
A Morse decomposition of 
a dissipative dynamical system can be obtained by dividing the phase space into grids 
and constructing a combinatorial multi-valued map over the grids \cite{Arai_etal_2009, Bush_etal_2012}.
In the case of time-series data generated by a dynamical system,
a combinatorial multi-valued map over the grids can be similarly constructed.
However,
time-series 
data obtained from real measurements (e.g., meteorological data)
are often highly stochastic due to the presence of noise.
A multi-valued map is then determined statistically by preferable transitions between the grids.
We consider time-series data taken from the first 
two principal components of the pressure patterns in the troposphere and the stratosphere 
in the northern hemisphere, measured over 31 years, 90 days in each 
year $\times$ every 6 hours per day.  The application of the MGSTD method to the
troposphere data
yields some particular transitions between the Morse sets,
corresponding to specific motions in the phase space spanned by the principal components.
The motions detected by our analysis are consistent with
changes between characteristic pressure patterns
that have been previously recognized in meteorological studies.
A similar result is also obtained with the stratosphere data.
\end{abstract}



\pagestyle{myheadings}
\thispagestyle{plain}

\section{Introduction}

\subsection{General background}
The study of dynamics based on time-series data obtained from
measurements of nonlinear phenomena has developed since the 1970s.
The seminal method of delay-coordinates by Ruelle and Packard et al.\cite{Packard_Crutchfield_Farmer_Shaw_1980}
led to the mathematical theory of reconstructing attractors from time-series data
by Aeyels \cite{Aeyels_1981}, Takens \cite{Takens_1981}, and Sauer, Yorke, and Casdagli 
\cite{Sauer_Yorke_Casdagli_1991}.
These outcomes have been successfully 
applied to obtain dynamical information of a great variety of nonlinear phenomena.

The present study examines
meteorological time-series data from the viewpoint of 
dynamical systems, and
proposes a new approach for extending 
such time-series analysis to detect not only attractors 
but also unstable dynamics,
by concatenating a set of time-series data
derived from scattered initial conditions in the phase space of a dynamical system.
The idea is based on a topological computation method
to obtain a so-called the Morse decomposition,
that is,
a decomposition of the phase space of the dynamics into finite numbers of 
isolated invariant sets, called Morse sets, that are related in a gradient-like manner.
This decomposition
may be considered as a crude but global representation of the entire dynamics in the phase space.
Recently developed computer-assisted methods for studying dynamics
\cite{AUTO, GAIO, CAPD} enabled us to understand various aspects of the dynamics of
concrete nonlinear systems.
Among these methods, 
several authors including one of the present authors have proposed 
a computational approach \cite{Arai_etal_2009, Bush_etal_2012}
for obtaining
Morse decompositions together with a concise description of the dynamics of 
each Morse set (in terms of the Conley index) of a given dynamical system.
This method was mainly developed for iterated maps with parameters,
but can also be applied to ordinary differential equations \cite{Miyaji_etal_2016}.

The basic idea of obtaining Morse decomposition of a given dynamical system
is to set a finite grid decomposition on a domain of interest in the phase space 
of the dynamics, and then to construct a combinatorial multi-valued map over 
the finite grid elements, which is a rigorous%
\footnote{
The multi-valued maps used in this paper
are constructed from time-series data
generated from meteorological measurements, and therefore
the results are not guaranteed to be mathematically rigorous.
}
outer approximation of the true time evolution of the dynamical system.  
Such a combinatorial multi-valued map may be represented as a finite directed graph,
with nodes representing the grid elements and edges representing the time evolution from 
one grid element to a set of grid elements that intersect with (outer-approximation of) 
the true image of the given grid element.  A strongly connected path component of 
this finite directed graph forms an isolated invariant set (or, more precisely, 
an isolating neighborhood whose maximal invariant subset is an isolated invariant set),
and the remaining part of the directed 
graph becomes gradient-like, since 
it contains no recurrent paths.

This idea is applicable not only to the numerically computed time evolution of dynamical systems,
but also to time-series data
obtained from measurements of phenomena driven by unknown dynamics,
if the set of time-series data is sufficiently large to
capture all the essential dynamical features in the phase space.
Provided this assumption holds,
one may then construct
a combinatorial multi-valued map from the time-series data in a similar way,
and hence obtain a Morse decomposition of the underlying dynamics directly
from the measurement,
rather than relying on mathematical models.

\subsection{Morse decomposition of the dynamics from stochastic time-series data 
and its application to meteorological time-series data}
The dynamics of real phenomena is inevitably subjected to noise.
Time-series data may thus be considered to contain some information
of time evolution governed by both a deterministic dynamical system and noise.
A reliable methodology is therefore essential
for understanding the underlying dynamics
and drawing meaningful conclusions.
This may involve removing the effects of noise to obtain
the deterministic contribution.

Two of the present authors and their collaborators are currently developing such 
a method \cite{Kokubu_Morita_Nomura_Obayashi},  which we call the {\em Morse 
graph method for stochastic time-series data (the MGSTD method)}.
Stochastic time-series in the phase space constitute an ensemble of transitions between grid elements.
Extracting statistically preferred and relevant transitions yields 
combinatorial multi-valued map,
for which Morse decomposition is performed
as described above.
When applied to time-series data generated 
from simple deterministic dynamical system models
with added stochastic terms,
the MGSTD method successfully reproduce the stable and unstable invariant sets,
as well as their connecting orbits,
of the noise-less deterministic dynamical systems.
A concise demonstration is given in Subsection \ref{sect:model} below.

A natural problem of interest is to apply the MGSTD method not only to the
time-series data generated from dynamical system models but also to those obtained
from measurement of real phenomena.
The purpose of this paper is to apply
the method to time-series data obtained from meteorological measurement.
The reason that, among other stochastic phenomena, we focus on the meteorological 
dynamics is as follows.
First, as evidenced by the success in weather forecasting,
basic physical equations governing meteorological behavior
are established,
such as the Navier-Stokes equation
in the Earth's rotational frame,
the continuity equation for dry air and other materials, and
the radiative transfer equation \cite{Holton_Hakim_2012}.
This enables us not only 
to predict a future state of meteorological variables to some extent, but also to 
diagnose the past and present states.  Second, artificial satellites' highly-frequent, 
spatially-dense, global observation
has provided reliable atmospheric data with 
sufficient spatial and temporal resolution since the 1980s.  By applying 
data assimilation techniques, physically consistent gridded data can be
created from observation and weather forecast data \cite{Kalney_2002}.
Third,
projecting the high-dimensional meteorological time-series data
onto a limited-dimensional phase space
gives globally stable orbits
with a non-Gaussian probability density function \cite{Kimoto_Ghil_1993a},
which implies that the motion is not completely random.
In fact,
although the existence of dynamically stable points was recently denied \cite{Stephenson_Hannachi_ONeill_2004},
the existence of preferable paths
between characteristic points corresponding to atmospheric states
has been indicated with the help of meteorological knowledge
\cite{Kimoto_Ghil_1993b,Luo_Cha_Feldstein}.
Some deterministic aspects can therefore be expected to exist in the stochastic time-series data.

\subsection{Outline of the paper}
The outline of this paper is as follows.
Section 2 introduces the MGSTD method.
First,
the mathematical theory of combinatorial Morse decomposition
of dynamical systems
and that of deterministic time-series data are reviewed.
Next, we review the method for constructing a multi-valued map from stochastic time-series data.
Examples of application to simple dynamical systems with noise
are briefly presented to demonstrate the usefulness of the MGSTD method.
The meteorological data, to which the MGSTD method is applied,
are then presented.
Section 3 reports the result of the MGSTD analysis of the datasets of the troposphere and the stratosphere.
Some specific motions in the phase space are observed
within highly stochastic time-series data.
We examine the dependence of the results
on the choice of parameter values,
to see the robustness of the results.
Section 4 is
devoted to discussion and concluding remarks.
The relevance of the observed motions
to existing meteorological knowledge is discussed,
and the results are shown to be consistent with, and complement to, earlier studies.

\section{Method}
\subsection{Morse decomposition of global dynamics from time-series data}
We here introduce the method for obtaining a Morse graph
from a given time-series dataset,
which is expected to represent the Morse decomposition
of the dynamical system
that underlies the measured time-series data.
We first briefly summarize the computer-assisted method,
developed in \cite{Arai_etal_2009, Bush_etal_2012},
for obtaining a Morse decomposition of a dynamical system
represented by a directed graph called a Morse graph.  
Section 2.1.4 applies this idea to time-series data
to obtain a Morse graph of the underlying dynamical system.

\subsubsection{Morse decomposition of a dynamical system}
As explained above, the Morse decomposition of a dynamical systems is a decomposition 
of the phase space into recurrent part and gradient-like part.  In this paper, we mainly 
consider discrete time dynamical systems, namely, an iterated map.  Let $X$ be a 
compact metric space and $f: X \to X$ a continuous map.

\begin{definition}[see \cite{Conley_1978}]\rm 
A {\em Morse decomposition} of the map $f$ is a finite collection of disjoint isolated invariant sets $S_1,\ldots,S_n$ (called {\em Morse sets}) with a strict partial ordering 
$\prec$ on the index set $\{1,\ldots,n\}$ such that,
for every $x \in X\setminus 
\cup_{i=1}^n S_i$ and every complete orbit $\gamma=\{x_n\}_{n\in{\mathbb Z}}$ 
through $x_0=x$,  i.e. $x_{n+1}=f(x_n)$ for all $n\in {\mathbb Z}$, there exist indices 
$i \prec j$ such that $x_n \to S_i$ and $x_{-n} \to S_j$ as $n \to +\infty$.  (In this case, 
$\gamma$ is called a {\em connecting orbit} from $S_j$ to $S_i$.)
\end{definition}

Here,
an {\em isolated invariant set} $S(\subset X)$ of $f$ means that it is an invariant set 
which has a compact neighborhood $N$ for which $S$ is its maximal invariant set in $N$ 
and sits in its interior, namely $S\subset \mathrm{Int}\; N$.  The neighborhood $N$ is 
called an {\em isolating neighborhood} of $S$.

Notice that Morse decomposition of a given dynamical system is not unique in general.  
The coarsest Morse decomposition of a map $f: X\to X$ consists of a single set $S$ 
which is the maximal invariant set of $f$ in $X$.  If $i$ and $j$
are indices such that $i \prec j$ 
but there is no other index $k$ such that $i \prec k \prec j$,
then a coarser Morse decomposition can be created
by replacing $S_i$ and $S_j$ with $S_i \cup S_j \cup 
C(i,j)$, where $C(i,j)$ denotes the union of all connecting orbits from $S_j$ to $S_i$.

For two Morse decompositions ${\mathcal S}=\{S_1,\ldots,S_n\}$ and ${\mathcal T}=
\{T_1,\ldots,T_m\}$, we say that ${\mathcal S}$ is a refinement of ${\mathcal T}$, if $n\geq m$ 
and if there is a surjective map $\iota:\{1,\dots,n\}\to\{1,\dots,m\}$ such that $S_i
\subset T_{\iota(i)}$ for any $i=1,\dots,n$.  By definition, any connecting orbit $\gamma$ 
between $S_i$ and $S_j$ is also contained in $T_k$ if $\iota(i)=\iota(j)=k$.

A Morse decomposition
${\mathcal S}=\{S_1,\ldots,S_n\}$ with a partial order $i\prec j$
can be represented in terms of a directed graph $G = (V,E)$,
where $V = \{S_1,\ldots,S_n\}$ and $(S_j,S_i) \in E$ iff $i \prec j$.  This graph is called 
a {\em Morse graph}.  In order to represent the computed Morse decomposition in 
a compact way, it is convenient to plot a Morse graph whose edges are determined by 
the {\em transitive reduction} of the relation $\prec$ which is a minimal relation $\prec'$ 
whose transitive closure retrieves $\prec$.
Such a representation is used below.

\subsubsection{Graph representation of dynamics}
To obtain a Morse decomposition of a map $f: X\to X$ with the aid of computer, 
we follow the idea of graph representation of the dynamics using a grid decomposition of the phase space, as given in \cite{Arai_etal_2009}.
In the case where $X$ is a compact 
domain in ${\mathbb R}^d$, 
and the map $f: X\to X$ is given by a mathematical formula using functions that can be 
handled by numerical computation,
let $\mathcal{Q}$ be a cubical grid decomposition of ${\mathbb R}^d$ that covers $X$,
and we aim to compute its image $f(Q)$ for $Q\in\mathcal{Q}$ by computer.
It is not in general possible to 
obtain the exact image by computer, and we can only expect to obtain its numerical 
approximation.
However,
the development of validated numerical computation has provided
various techniques for numerically computing rigorous outer bounds.
If one can compute a rigorous outer approximation 
$[f(Q)]$ of the image of a grid element $Q\in {\mathcal Q}$ under $f$
(e.g., by using interval arithmetics),
let ${\mathcal F}(Q)$ be the set of all grid elements in ${\mathcal Q}$
that intersect with $[f(Q)]$, 
namely ${\mathcal F}(Q) = \{ R\in \mathcal{Q}\mid R\cap [f(Q)]\neq\emptyset\}$.  
This defines a  multi-valued map ${\mathcal F}$ from ${\mathcal Q}$ to itself.  By definition, 
the union of all grid elements in ${\mathcal F}(Q)$ completely contains the true image 
$f(Q)$.  In this sense,  ${\mathcal F}$ can be considered as a rigorous outer-approximation 
of the map $f:X\to X$.  We call ${\mathcal F}$ a {\em combinatorial representation} of $f$.  
Since it is a multi-valued map on ${\mathcal Q}$,  we use the notation $\mathcal{F}: 
\mathcal{Q} \multimap \mathcal{Q}$ to distinguish it from a usual single-valued map.

A combinatorial representation $\mathcal{F}: \mathcal{Q} \multimap \mathcal{Q}$ of $f$ 
can be equivalently represented by means of a directed graph $G = (V,E)$,
where
$V = \mathcal{Q}$ and $(Q,Q') \in E$ iff $Q' \in \mathcal{F} (Q)$.  The analysis of $G$ 
provides information on the asymptotic dynamics of $f$ represented by $\mathcal{F}$.
For instance, each {\em combinatorial invariant set}
defined as a set $\mathcal{S} \subset 
\mathcal{Q}$ for which $\mathcal{S} \subset \mathcal{F} (\mathcal{S}) \cap 
\mathcal{F}^{-1} (\mathcal{S})$
represents an isolating neighborhood $|\mathcal{S}|$ with respect to $f$,
where $|\mathcal{S}|$ stands for the geometric realization of the set $\mathcal{S}$ of 
the collection of grid elements. 
Moreover, a {\em combinatorial attractor}
defined as a set $\mathcal{A} \subset \mathcal{Q}$
such that $\mathcal{F} (\mathcal{A}) \subset \mathcal{A}$
represents an isolating neighborhood $|\mathcal{A}|$ whose invariant part $A$ is stable 
in the sense of Conley \cite{Conley_1978}:
Every forward orbit starting from a point $x$ in some open neighborhood of $A$
(or, more precisely, in $\mathrm{Int} |\mathcal{A}|$) approaches 
$A$ ($\mathrm{dist} (f^n(x), A) \to 0$ as $n \to \infty$).
In particular,
if there exist two disjoint combinatorial attractors for $\mathcal{F}$,
then this implies the existence of two disjoint basins of attraction for $f$,
thus the dynamical system $f$ is (at least) bistable.

\subsubsection{Combinatorial Morse decompositions}
An extensive analysis of the dynamics with its combinatorial representation
can be performed by computing the {\em strongly connected path components}, which are defined 
in terms of a directed graph as the equivalent formulation of the combinatorial multi-valued map.
The {\em strongly connected path components} of a directed graph $G=(V,E)$ are the 
maximal sets of vertices $C \subset V$ 
that satisfy the following property:
for each $v, w \in C$, there exists a path from $v$ to $w$ with vertices in $C$ and also 
a path in the opposite direction (from~$w$ to~$v$) through $C$.  In 
\cite{Kalies_Mischaikow_Vandervorst_2005} it is shown that all the strongly connected path
components of $G$ form isolating neighborhoods for the union of all the chain recurrent 
sets of the dynamical system, and thus can serve as a {\em combinatorial Morse 
decomposition} $\{\mathcal{M}_i \mid i = 1, \ldots, k\}$, for some $k > 0$,
representing
a family of isolating neighborhoods $|\mathcal{M}_i|$.  The sets $\mathcal{M}_i$ are called 
{\em combinatorial Morse sets}.  A partial order $\prec$ between the computed 
combinatorial Morse sets can be determined by the analysis of paths in $G$ connecting 
those sets:  $i \prec j$ if $i \not = j$ and if there exists a path in $G$ from any vertex in 
$\mathcal{M}_j$ to any vertex in $\mathcal{M}_i$.  With the use of the graph $G$, 
a combinatorial Morse decomposition can be computed by algorithms introduced in 
\cite{Arai_etal_2009,Ban_Kalies_2006,Bush_etal_2012}.
In the following,
we refer to {\em Morse sets} as combinatorial Morse sets $\mathcal{M}_i\in\mathcal{Q}$ and 
their geometrical representation $|\mathcal{M}_i|\in\mathbb{R}^d$ without distinction.

\subsubsection{Combinatorial Morse decompositions from deterministic time-series data}
Supposing a time-series dataset is derived from an unknown dynamical system,
we seek to obtain a Morse decomposition of the underlying dynamical system.
More precisely, let $f:X\to X$ be a continuous 
map of a compact domain $X\subset {\mathbb R}^d$, and let $\pi:X\to {\mathbb R}^m$ 
be an observation map.
Then, for each initial point $x\in X$,
we obtain a finite time-series 
$\{y_n\}_{n=0,\dots,N}$
given by
$y_n=\pi(f^n(x))$.
For a finite set
$\Xi=\{x^k\}_{k\in K}$
of initial conditions in $X$,
where $K$ is a finite index set,
and a set of natural numbers
$\{N_k\}_{k\in K}$,
we obtain 
a set $D$ of finitely many finite time-series data as follows:
\begin{align}
\label{eqn:data_set}
D=\{y_n^k=\pi(f^n(x^k))\in{\mathbb R}^m  \mid k\in K, \ n=0, 1,\dots,N_k \}
\end{align}
This corresponds to the image of points
$\tilde\Xi=\{f^n(x^k)\}_{k\in K, \ n=0,\dots, N_k}$ 
under the observation function $\pi$.
For later purpose, we also define 
the subset
\begin{align}
\label{eqn:data_subset}
D'=\{y_n^k=\pi(f^{n}(x^k))\in{\mathbb R}^m \mid k\in K, \ n=0, 1,\dots,N_k-1 \}
\end{align}
of $D$.

Let $\mathcal{R}$ be a cubical grid decomposition of a domain in ${\mathbb R}^m$ that 
covers $\pi(X)$, and define $\mathcal{Y}=\{R\in\mathcal{R}\mid R\cap D\not=\emptyset\}$, 
and $\mathcal{Y}'=\{R\in\mathcal{R}\mid R\cap D'\not=\emptyset\}$.  Then we define 
a combinatorial multi-valued map $\mathcal{F}: \mathcal{Y}' \multimap \mathcal{Y}$ as follows:  for $R\in{\mathcal{Y}}$ and $R'\in{\mathcal{Y}'}$, we define $R\in{\mathcal F}(R')$ 
iff there exist
$y^k_n, y^k_{n+1}\in D$
such that
$y^k_n\in R'$ and $y^k_{n+1}\in R$
hold.
Note that the multi-valued map $\mathcal{F}$ should be regarded as a self-map $\mathcal{F}: 
\mathcal{Y} \multimap \mathcal{Y}$ as in \S 2.1.2.
However, in practice, its domain of 
definition is restricted to $\mathcal{Y}'\subset \mathcal{Y}$ due to the finiteness of the 
time-series data.

Once $\mathcal{F}$ is given, the same procedure as described
above can be followed,
yielding a finite collection
$\{\mathcal{M}_p \mid p = 1, \ldots, P \}$
of the strongly connected path components of the equivalent directed graph 
representations of $\mathcal{F}$, and hence the Morse graph of $\mathcal{F}$.

We then take a pullback of the sets
$\{|\mathcal{M}_p| \mid  p = 1, \ldots, P \}$,
namely, define
$N_p=\pi^{-1}(|\mathcal{M}_p|)$.
Obviously, these are disjoint compact subsets of 
$X$, 
since the space $X$ is assumed to be compact.
In the present situation, we assume that
the unknown dynamical system is dissipative,
and hence that it has a so-called global attractor,
which is the maximal compact invariant set that attracts all the trajectories.
It is therefore reasonable in practice
to assume that $X$ is compact, even in the case of dynamical 
systems driven by (dissipative) PDEs.
Provided the time-series data are sufficiently abundant,
one can show that these sets are indeed isolating neighborhoods
and that their maximal invariant sets $M_p = \mathrm{Inv}(N_p)$
give a Morse decomposition of the unknown 
dynamical system $f:X\to X$.
In this paper, we do not intend to give a detailed statement of the theory
nor the precise conditions for the theory to be applied,
mainly because these are unverifiable in practical applications.
We therefore simply note that there is some mathematical
basis for this statement, even if it is not always practically meaningful for application problems \cite{Kokubu_Morita_Nomura_Obayashi}.

\subsection{Multi-valued map from stochastic time-series data}
We here review a method for constructing a combinatorial multi-valued map
$\mathcal{F}: \mathcal{Y}' \multimap \mathcal{Y}$
when the time-series data is noisy,
as in the case of meteorological data
\cite{Kokubu_Morita_Nomura_Obayashi}.

Assume we are given a set $D$ of $m$-dimensional time-series data,
as in \S 2.1.4.
A sufficiently large $\tilde\Omega = [-\tilde L, \tilde L]^m$
is chosen to contain $D$.
Typically $\tilde L=4$ suffices for
the PCA data considered in this paper, as the PCA scores
are each normalized to have zero mean and unit standard deviation.
We put the grid decomposition $\mathcal{R}$
on the domain $\tilde\Omega$ with grid size $h>0$,
from which the combinatorial multi-valued map is constructed.
The map may depend on the setting of the grid decomposition,
e.g. the size and location of the grid elements.
In this paper,
we only consider a square grid decomposition $\mathcal{R}$ of the size $h$,
and we control its location by the shift parameters $\delta_\ell\in[0,h)\ (\ell=1,\dots,m)$.
Thus, we consider the grid decomposition over the domain 
$$
\Omega = [-L+\delta_1,L+\delta_1]\times[-L+\delta_2,L+\delta_2]\times\cdots\times[-L+\delta_m,L+\delta_m],
$$
where $L(<\tilde L)$ is a positive integer multiple of $h$ such that $D\subset\Omega$.  
Hence, a grid element in $\mathcal{R}$ takes of the form $Q^i=[a_1^i, a_1^i+h]\times
[a_2^i,a_2^i+h]\times\cdots\times[a_m^i,a_m^i+h]$ where $a_\ell^i-\delta_\ell$ is a positive 
integer multiple of $h$.

We determine a multi-valued map from the time-series as follows.
We first define the initial multi-valued map $\tilde{\mathcal{F}}$ with respect to
$\mathcal{R}$ from the dataset $D$ just as in \S 2.1.4, and then we modify 
$\tilde{\mathcal{F}}$ so as to reflect the effect of noise by taking into account of 
several characteristic numbers introduced below.  Let $\nu_i$ be the number of 
data points in a grid element $Q^i\in\mathcal{R}$, namely $\nu_i=\#(Q^i\cap D)$, and 
$\mu_{i\to j}$ the number of transitions from $Q^i$ to $Q^j$, namely
$\mu_{i\to j}=\#\{(y^k_n,y^k_{n+1}) \mid y^k_n\in Q^i\cap D', \ y^k_{n+1}\in Q^j\cap D\}$.
We then define the transition probability, or the conditional probability, 
from $Q^i$ to $Q^j$ by
\begin{align}
T_{i\to j}=\mu_{i\to j}/\nu_i
\label{eqn:Def_Tij}
\end{align}
We consider the transition between grid elements $Q^i$ and $Q^j$
to be determined by
the transition probability, and by a parameter $\rho(\ge 1)$ indicating the degree of 
superiority, as follows.
If $\rho < T_{i\to j}/ T_{j\to i}$,
then we say the transition from $Q^i$ to $Q^j$ is superior to
the opposite transition from $Q^j$ to $Q^i$, denoted by $i\to j$.
Similarly,
we say $Q^j$ to $Q^i$ is superior to the opposite, denoted by $j\to i$,
if $\rho < T_{j\to i}/ T_{i\to j}$, or equivalently, $T_{i\to j}/ T_{j\to i} < \rho^{-1}$.  Otherwise, 
we say the transitions between $Q^i$ and $Q^j$ are comparable, denoted by
$i\leftrightarrow j$.
In summary,
\begin{subequations}
\begin{alignat}{4}
i\to j &\quad\mbox{if}\quad\;\;\; \rho &<  T_{i\to j}/ T_{j\to i} &
\\
i\leftrightarrow j &\quad\mbox{if}\quad \rho^{-1} &\le T_{i\to j}/ T_{j\to i} & \le \rho
\\
i\leftarrow j &\quad\mbox{if}\quad & T_{i\to j}/ T_{j\to i} &< \rho^{-1}
\end{alignat}
\label{eqn:Tij_rho}
\end{subequations}
Note that self transitions $i\to i$
are always taken into account unless $\mu_{i\to i}=0$.

\begin{figure}[htbp]
\centering
\tabcolsep=0mm
\begin{tabular}{ccc}
\begin{minipage}[c]{0.33\textwidth}
\includegraphics[width=\textwidth]{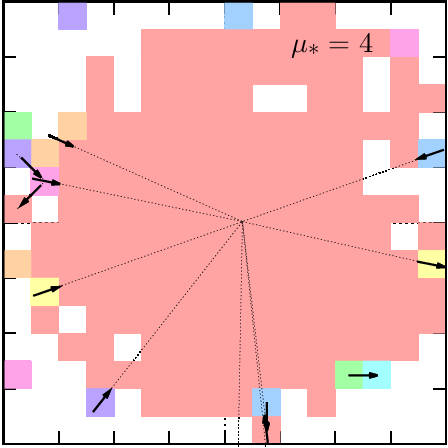}
\end{minipage}
&
\begin{minipage}[c]{0.33\textwidth}
\includegraphics[width=\textwidth]{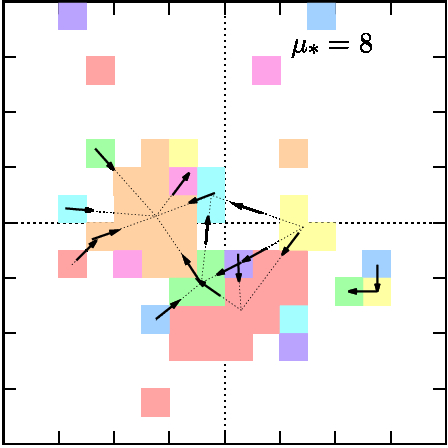}
\end{minipage}
&
\begin{minipage}[c]{0.33\textwidth}
\includegraphics[width=\textwidth]{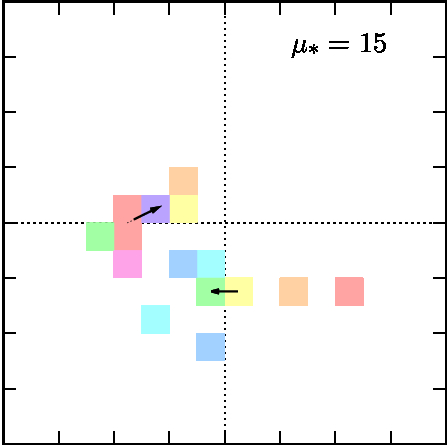}
\end{minipage}
\end{tabular}
\caption{
Example of showing the splitting of Morse sets with increasing $\mu_*$ in the troposphere,
with $h=0.25$, $\rho=1.1$, and $(\delta_1, \delta_2)=(0,0)$.
The abscissa and ordinate are PC1 and PC2, respectively,
in all panels.
}
\label{fig:MSvsmu}
\end{figure}

In addition, we avoid overestimating rare events; otherwise, for example,
the transition probability for $\mu_{i\to j}=10$ and $\nu_i=20$
would be regarded to be equal to that for $\mu_{k\to l}=1$ and $\nu_k=2$,
though the latter may occur just by chance.
To this end, we introduce a threshold $\mu_*$ such that only
the transitions $i\to j$ with $\mu_{i\to j}\ge \mu_*$ are taken into account;
this threshold is also applied to self transitions.
In fact, as shown below in \reffig{fig:MSvsmu},
when $\mu_*$ is small, almost all the grids are strongly path connected,
leading to a Morse set so large as to cover all the relevant domain in the phase space;
with increasing $\mu_*$, such a Morse set is divided 
into several Morse sets, some of which are robustly seen for various choices of the other 
parameters; when $\mu_*$ is too large, even relevant transitions are removed, leaving only 
independent Morse sets with no transitions between them.

\begin{figure}[htbp]
\tabcolsep=0.01\textwidth
\begin{tabular}{cc}
(a) Morse graph & (b) Morse sets
\\
\begin{minipage}{0.45\textwidth}
\includegraphics[width=\textwidth]{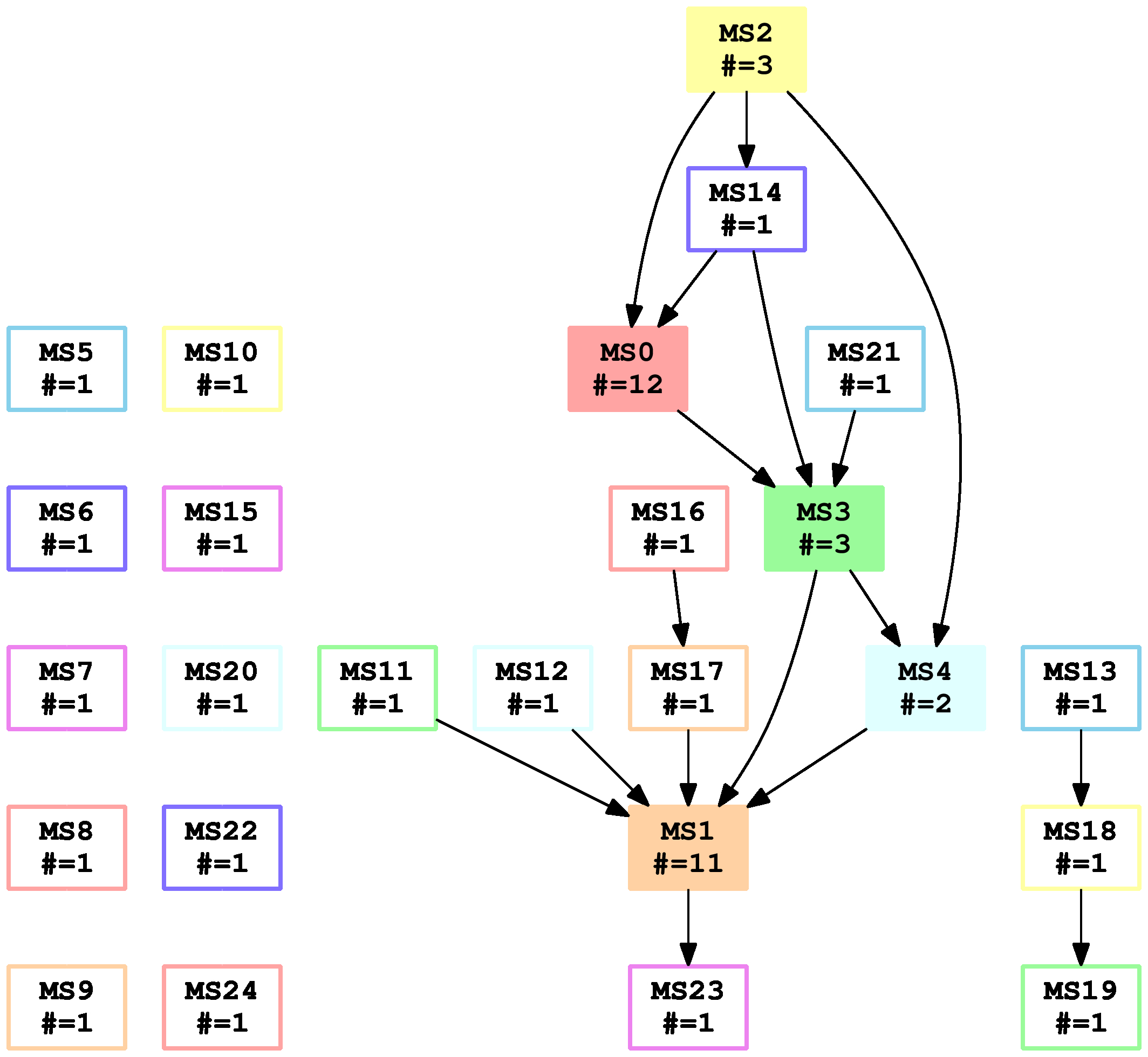}
\end{minipage}
&
\begin{minipage}{0.5\textwidth}
\includegraphics[width=\textwidth]{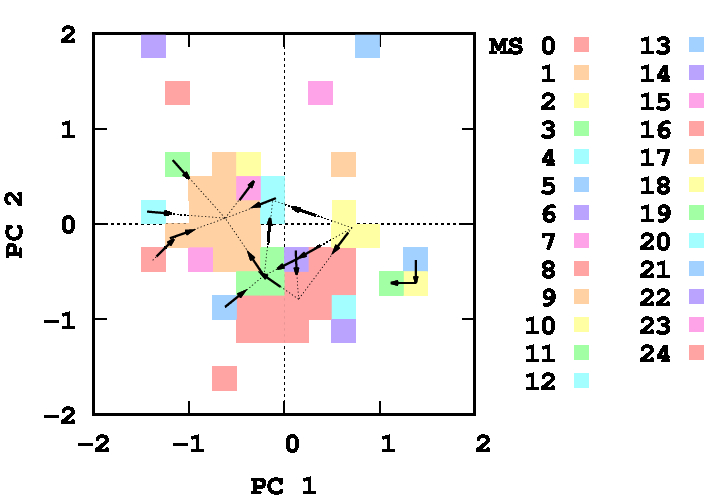}
\end{minipage}
\end{tabular}
\caption{
(a) Sample Morse graph and (b) the corresponding Morse sets
for the troposphere,
with $h=0.25$, $\rho=1.1$, $(\delta_1, \delta_2)=(0,0)$, and $\mu_*=8$.
Morse sets are denoted by {\rm MS}$i$ ($i=0,1,2,\dots$),
and ordered according to size, which also determines the color coding of the Morse sets.
Note that isolated Morse sets of unit size may be spurious
resulting from slow dynamics.
The color coding is the same in (a) and (b).
In (a), the Morse sets of size two or more are filled with color.
In (b), the arrows and dotted lines,
connecting the barycenters of Morse sets,
represent the gradient-like connections in the Morse graph in (a).
See also the MGSTD algorithm below.
}
\label{fig:Morse_graph_set_tropo}
\end{figure}

We thus define the multi-valued map $\mathcal{F}$ by 
$Q^j\in\mathcal{F}(Q^i)$, if $Q^j\in\tilde{\mathcal{F}}(Q^i)$, either $i\to j$ or 
$i\leftrightarrow j$, and $\mu_{i\to j}\ge \mu_*$ are all satisfied.
\reffig{fig:Morse_graph_set_tropo} is a sample Morse graph and the corresponding
Morse sets, in the case of $m=2$,
obtained from the meteorological time-series data explained below in \S 2.4.
Notice that the Morse graph in \reffig{fig:Morse_graph_set_tropo}(a) can be recovered using
the arrows in \reffig{fig:Morse_graph_set_tropo}(b), which also exhibit the location of
the corresponding Morse sets in the phase space.
We call such a presentation of the Morse graph 
the {\em phase space presentation of the Morse graph},
and we adopt it in the following.

We thus have 
$4+m$ 
control parameters: $m, h, \rho, \mu_*, \mbox{ and } \delta_1,\delta_2,\cdots,\delta_m$
for defining the multi-valued map $\mathcal{F}$.
Among them, we set the dimension of the \lq phase space' $m=2$ throughout the paper,
following earlier studies
\cite{Kimoto_Ghil_1993a,Kimoto_Ghil_1993b,Inatsu_Nakano_Mukougawa_2013,Inatsu_Nakano_Kusuoka_Mukougawa_2015}.

Two issues must be addressed to
extract information about the `deterministic' transitions arising from the unknown underlying dynamics
and distinguished from the stochastic time-series data:
(A) how to select parameters of computation, and (B) how to 
exhibit the aspects of transitions from the results of computation.
Our approach, explained below, is called
the {\em Morse graph method for stochastic time-series data} 
(abbrev. as the {\em MGSTD method}).

With regard to the issue (A), after fixing
the number of principal components $m=2$,
we vary the other $m+3=5$ parameters,
namely the grid size $h$, the degree of superiority
$\rho$ for determining the direction of transition 
between grid elements, the threshold $\mu_*$ for the number of transitions between grid 
elements, and $\delta_1,\delta_2\ (0\leq\delta_\ell<h, \ell=1,2)$, the position of the center 
of the grid.

An essential parameter is $\mu_*$.
Since the original data is highly stochastic,
a single large Morse set may be obtained as a result of stochastic recurrence in the data,
if $\mu_*$ is inappropriately set (See \reffig{fig:MSvsmu}).
For each pair of grid elements,
the transitions occurring less frequently than $\mu_*$ times are discarded,
considered to result from stochasticity.
As $\mu_*$ is increased,
fewer pairs of grid elements display transitions,
meaning a decrease in the chance of recurrence.
Upon reaching some value of $\mu_*$,
a single large Morse set splits into several smaller Morse sets of comparable size.

\begin{figure}[htbp]
\begin{minipage}[c]{0.49\textwidth}
\centering
(a) Troposphere
\includegraphics[width=\textwidth]{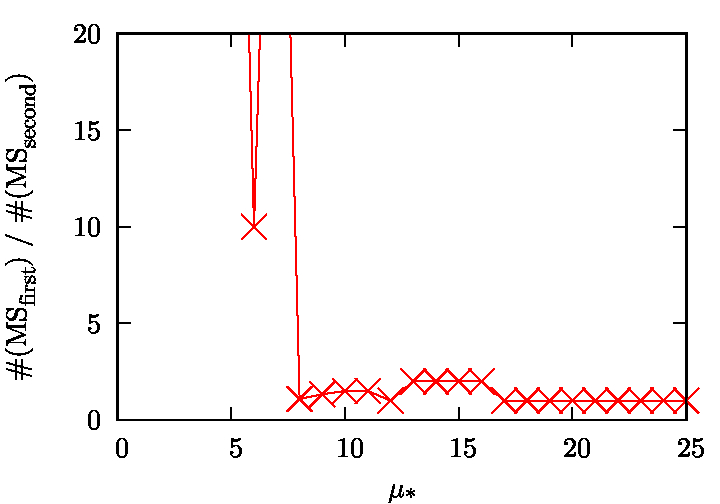}
\end{minipage}
\begin{minipage}[c]{0.49\textwidth}
\centering
(b) Stratosphere
\includegraphics[width=\textwidth]{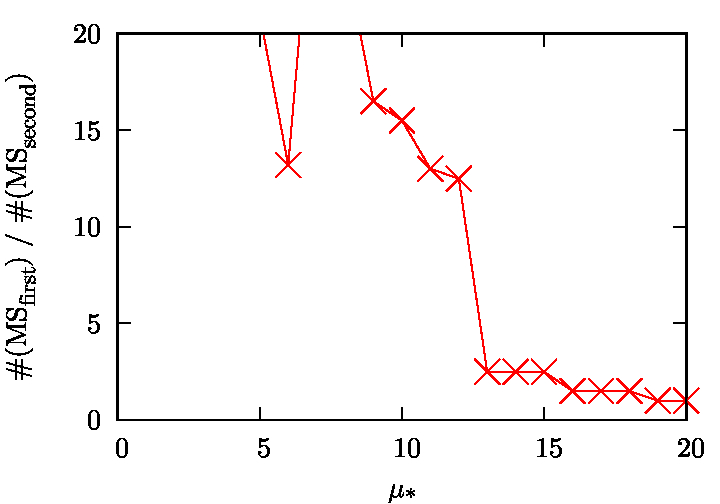}
\end{minipage}
\caption{
Dependence on $\mu_*$ of the ratio
$\#({\rm MS}_{\rm first})/ \#({\rm MS}_{\rm second})$
for (a) the troposphere and (b) the stratosphere;
$h=0.25$, $\rho=1.1$, $(\delta_1, \delta_2)=(0,0)$.
}
\label{fig:Rvsmu}
\end{figure}

It is therefore reasonable to set the smallest value of $\mu_*$
for which the large single Morse set splits into
several Morse sets of comparable size.
The value may be estimated by the ratio of the sizes of the two largest Morse sets.
More specifically, we select $\mu_*^{\circ}$ as
$$
\mu_*^{\circ}=
\min\left\{\mu_*\mid \frac{\#({\rm MS}_{\rm first})}{\#({\rm MS}_{\rm second})}<A\right\}
$$
where $\#({\rm MS}_{\rm first [resp. second]})$ is the size (i.e. the number of grid elements) of
the first and second Morse sets (ordered according to size) for a given $\mu_*$,
and $A$ is a value suitably chosen from the data.

The above ratio is computed by varying $\mu_*$
for fixed values of $\rho$, $h$, and $\delta_\ell$.
A sharp drop is observed at some values of $\mu_*$ (\reffig{fig:Rvsmu}).
Using the above criterion,
the first value of $\mu_*$ after the sharp drop is selected,
giving $\mu_*^{\circ}=8$ for the troposphere and $13$ for the stratosphere.
In both cases, any value of $A$ between 3 and 9 gives
the same values for $\mu_*^{\circ}$.
We have therefore chosen $A=5$ for the computations throughout the paper.

All the above computations were done with $\rho=1.1$ and $h=0.25$;
see below for details of how to set these values.

\medskip

With regard to the issue (B),
we introduce the idea of a `vector field'
to display temporal transitions between Morse sets,
by the following algorithm.

\medskip
\noindent{\bf MGSTD algorithm}:
\begin{description}
\item{Step 0:} Fix $\mu_*=\mu_*^{\circ}$, and vary $\delta_1, \delta_2$ each from $0$ to $h$ with some increment (e.g., $0.01$).
\item{Step 1:} Compute the Morse graph and its Morse sets.
\item{Step 2:} For each gradient-like connection between Morse sets in the Morse graph, 
namely for each arrow from one Morse set ${\rm MS}_i$ to one of its descendent Morse
sets, say ${\rm MS}_j$,
compute the vector $v_{ij}$
that is parallel to the vector $p_i-p_j$,
where $p_k\ (k=i,j)$ is the barycenter of $|{\rm MS}_k|$,
and such that its length $|v_{ij}|=(\#({\rm MS}_i)+\#({\rm MS}_j)/2$.
The center of  $v_{ij}$ is placed
at the $\#({\rm MS}_i) : \#({\rm MS}_j)$ interpolation point $q_{ij}$
of $p_i$ and $p_j$.
We thus obtain the distribution of vectors $v_{ij}$ for all the arrows in the Morse graph.
\item{Step 3:} For a choice of $\delta_1, \delta_2$, and a grid element $R$,
take an `average' of the vectors $v_{ij}$
whose centers belong to $R$,
by dividing their sum by their number and placing it at the center of $R$.
This vector is denoted by $v_R(\delta_1, \delta_2)$.
This gives a  distribution of the vectors over the grid $\mathcal{R}$,
for $\delta_1, \delta_2$.
\item{Step 4:} Take an `average' of the distribution of the vectors
by varying $\delta_1, \delta_2$ as follows.
For each grid element $Q$ in the canonical grid decomposition $\cal{Q}$ with $\delta_1=0, \delta_2=0$,
compute the average $w(Q)$ of all the vectors $v_R(\delta_1, \delta_2)$
whose centers belong to $Q$.
Place the vector $w(Q)$ at the center of the grid element $Q$.  
\end{description}
\medskip

We thus obtain the `vector field' $\{w(Q)\mid Q\in\cal{Q}\}$
over the canonical grid decomposition $\cal{Q}$,
which we call the {\em MGSTD vector field}.  In \S 2.3, we shall 
apply the MGSTD method to some mathematical models and compute the Morse graphs 
as well as the MGSTD vector fields.

\bigskip

\begin{figure}[bhtp]
\begin{minipage}[c]{0.49\textwidth}
\centering
(a) Troposphere
\includegraphics[width=\textwidth]{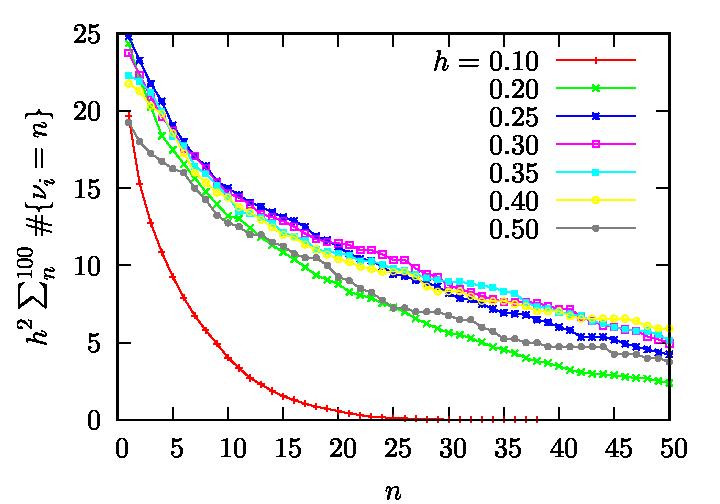}
\end{minipage}
\begin{minipage}[c]{0.49\textwidth}
\centering
(b) Stratosphere
\includegraphics[width=\textwidth]{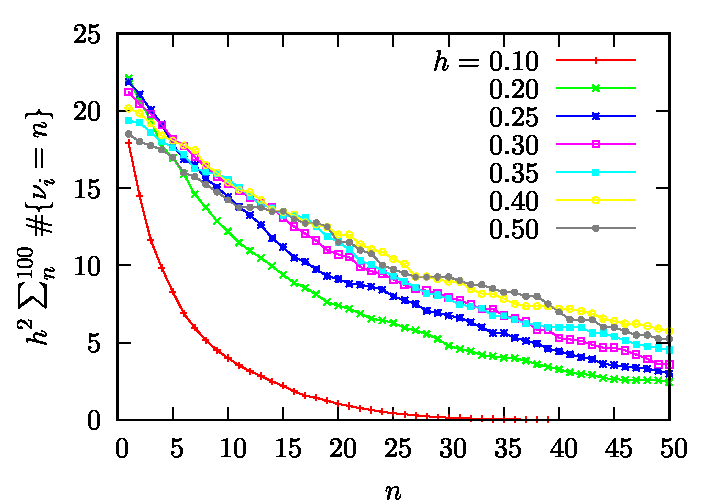}
\end{minipage}
\caption{
The relative portion of the grid elements of the size $h$ in the domain of the 
PC1-PC2 plane that carry $n$ or more data points,
for (a) the troposphere and (b) the stratosphere.
}
\label{fig:hselect}
\end{figure}

The grid size $h$ suitable for extracting the important aspects of dynamics
hidden in the time-series data is determined as follows.
Since the unity of PC1 and PC2 is their standard deviation
in which most of the data points are,
each unit square of the PC1-PC2 plane should be divided at least several times,
e.g., into $3^2$, $4^2$, or $5^2$ pieces,
by choosing small values of $h$, e.g., $h=1/3, 1/4$ or $1/5$, respectively.
Most of the grid elements, on the other hand, should contain sufficiently many number of 
data points, so that the transition probability among grid elements is relevant.  Smaller value 
of $h$ implies, however, smaller number of data points in a grid element.
A suitable lower bound for $h$ that meets the above two requirements
is determined by considering
$h^2\sum_n^{100}\#\{\nu_i=n\}$,
i.e., the relative portion of the grid elements of size $h$ in the domain of the PC1-PC2 plane that 
carry $n$ or more data points.  
Here, the number $n=100$ in the summation is practically considered as infinity.  
The dependence of this quantity on $n$ is shown in \reffig{fig:hselect}
for varying $h$.
By regarding between $10$ and $20$ data points
as sufficient in a grid element, the relative portion
is maximal at around $h=0.25$ for the troposphere,
and at around $h=0.3$ for the stratosphere.
This meets the above requirements, namely the unit square of the PC1-PC2 plane is divided 
into $4^2$ and $3^2$ grid elements for the troposphere and the stratosphere, respectively,
and the relative ratio of the grid elements
that contain between $10$ and $20$ data points is maximal. 
In order to see the robustness of the computation results to the value of 
$h$, we also vary $h$ in the range of $0.2\leq h\leq 0.3$ for the troposphere, and of 
$0.25\leq h\leq 0.35$ for the stratosphere.

The remaining parameter to be determined is $\rho$.
The value should not be very different from the unity, 
otherwise unnecessarily many pairs of grid elements would be regarded to have bidirectional transitions,
hence preventing the detection of meaningful dynamics in the data.
We therefore adopt $\rho=1.1$, $1.3$, and $1.5$ for the computations.

\subsection{Application to mathematical models}
\label{sect:model}
We here briefly demonstrate the application of the MGSTD method
to simple dynamical systems with noise.
Further details of the analysis will appear elsewhere
\cite{Kokubu_Morita_Nomura_Obayashi}.

\begin{figure}[htbp]
\centering
\parbox[b][10ex][t]{2em}{(a)}
\includegraphics[width=0.45\textwidth]{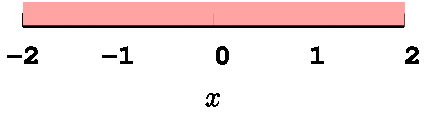}
\\
\vspace{\baselineskip}
\parbox[b][10ex][t]{2em}{(b)}
\includegraphics[width=0.45\textwidth]{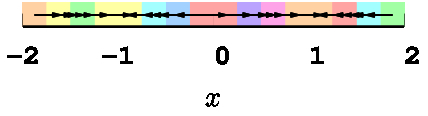}
\\
\vspace{\baselineskip}
\parbox[b][10ex][t]{2em}{(c)}
\includegraphics[width=0.45\textwidth]{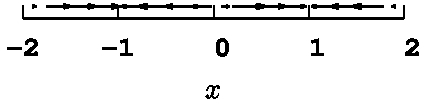}
\caption{
(a) Morse set for \refeqn{eqn:sde} with the dataset
$D_1^{\rm 1D}$
without parameter 
setting other than $h=0.25$,  (b) Morse sets and the phase space presentation of the 
Morse graph for \refeqn{eqn:sde}  with the dataset
$D_1^{\rm 1D}$
with
$h=0.25$, $\rho=1.1$, $\delta_1=0$, and $\mu_*^\circ=2$,
(c) MGSTD vector field 
$\{w(Q)\}$ for \refeqn{eqn:sde} with the dataset
$D_2^{\rm 1D}$
with $h=0.25$, $\rho=1.1$,
and
$\mu_*^\circ=1$ (the average is $1$).
}
\label{fig:dw}
\end{figure}

We consider the one-dimensional stochastic differential equation,
\begin{align}\label{eqn:sde}
dx_t&=x_t(1-x_t^2)dt+\sigma dB_t,
\end{align}
where
$B_t$ denotes a standard one-dimensional Brownian motion \cite{Oksendal}.
Without noise ($\sigma=0$), there exist two stable fixed points at $x=\pm 1$
and one unstable fixed point at $x=0$.
We set $\sigma=\sqrt{0.2}$ in the following.

We first consider the case where sufficiently many number of data are given;
we take $1000000$ initial points
uniformly randomly distributed over the interval $[-2,2]$,
and, for each initial point,
we make one step time evolution
by integrating \refeqn{eqn:sde}
over the time interval $[0,\Delta t]$
with a time step of $\Delta t=0.1$.
We performed the integration using the stochastic Runge-Kutta method
\cite{Honeycutt1992} with a time increment 0.001.
Thus we obtain the dataset of the form \refeqn{eqn:data_set},
$$
D^{\rm 1D}_1=\{ \xi_n^k = x^k_{n\Delta t} \mid k=1,\dots,1000000, \ n=0,1\}.
$$
When the multi-valued map is determined without using the parameter setting given in \S 2.2,
the resulting Morse graph becomes as shown in \reffig{fig:dw}(a),
with one large Morse set covering the relevant region of the one-dimensional phase space.
When the multi-valued map is determined by applying the MGSTD method
with the above criterion of parameter setting,
the Morse graph changes to that shown in \reffig{fig:dw}(b),
with two Morse sets at around $x=1$ (orange) and $x=-1$ (yellow),
corresponding to the sink vertices of the Morse graph,
and a Morse set at around $x=0$ (red),
corresponding to the source vertex of the Morse graph.

We next consider the case where insufficiently many number of data are given,
just like the case of the meteorological data in \S 2.4; 
we take $30$ initial points uniformly randomly distributed over the interval $[-2,2]$,
and make $399$ steps time evolution by integrating \refeqn{eqn:sde}
over the time interval $[0, 10-\Delta t]$ with a time step of $\Delta t=0.025$.
These numbers are almost as many in number as
the number of the meteorological data sets that we are going to analyze.
The integration scheme is the same as the above.
The dataset thus obtained is,
$$
\Tilde{D}_2^{\rm 1D} = \{ x^a_t \in\mathbb{R} \mid
a=1,\dots, 30, \ t=n\Delta t, \ n=0, 1, \dots, 399 \}.
$$
We consider this as the dataset of the form \refeqn{eqn:data_set},
$$
D_2^{\rm 1D}=\{ \xi_n^k \in\mathbb{R} \mid k\in K, \ n=0,1,\dots,99\}
$$
with the index set $K = \{ (a,c) \mid a=1,\dots,30, \ c=0,1,2,3\}$,
by $\xi^k_n=x^a_t$ with $t=0.1n+0.025c$.
Due to the shortage of the number of data points, the Morse graph in  this case
depends strongly on the choice of the origin of grid $\delta_1$.
Even in such a case, however, the MGSTD vector field $\{w(Q)\}$ shows a similar structure of the 
aspects of transitions between Morse sets as shown in \reffig{fig:dw}(c), namely
there are accumulation of arrows around $x=1$ and $-1$, and separation at around 
$x=0$.
We can therefore say that the MGSTD method qualitatively reproduces the deterministic nature of \refeqn{eqn:sde} without noise.

\begin{figure}[htbp]
\begin{flushleft}
\parbox[b][35ex][t]{0.04\textwidth}{(a)}
\includegraphics[width=0.44\textwidth]{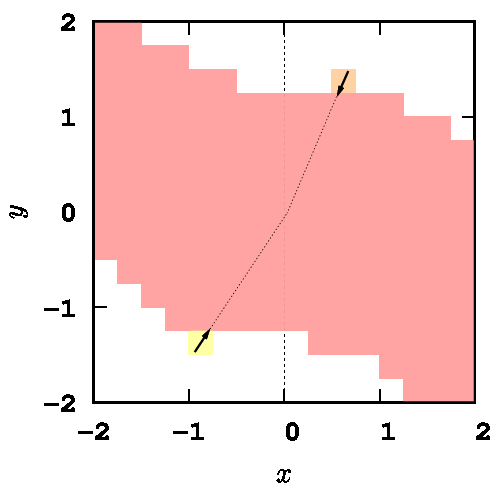}
$\;$
\parbox[b][35ex][t]{0.04\textwidth}{(b)}
\includegraphics[width=0.44\textwidth]{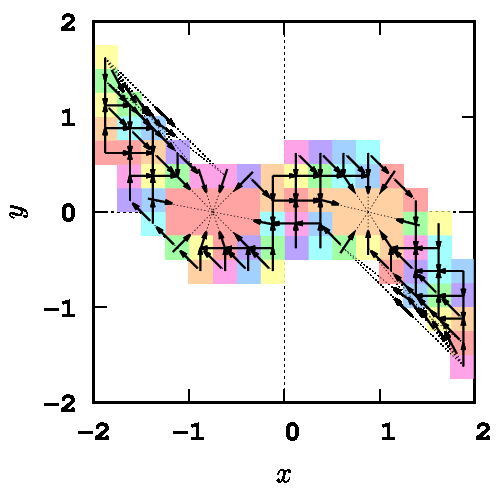}
\\
\parbox[b][35ex][t]{0.04\textwidth}{(c)}
\includegraphics[width=0.44\textwidth]{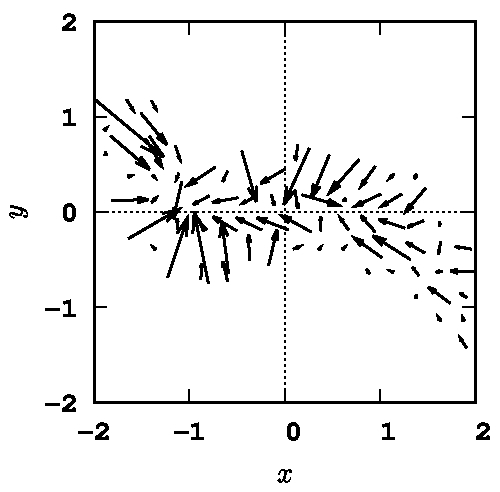}
\end{flushleft}
\caption{
(a) Morse sets and the phase space presentation of the Morse graph for 
\refeqn{eqn:sde2}  with the dataset
$D_1^{\rm 2D}$
without parameter setting other than $h=0.25$,
(b) Morse sets and the phase space presentation of the Morse graph for 
\refeqn{eqn:sde2} with the dataset
$D_1^{\rm 2D}$
with $h=0.25$, $\rho=1.1$, 
$(\delta_1, \delta_2)=(0,0)$, and $\mu_*^\circ=461$,
(c) MGSTD vector field $\{w(Q)\}$ 
for \refeqn{eqn:sde2} with the dataset
$D_2^{\rm 2D}$
with $h=0.25$, $\rho=1.1$, and $\mu_*^\circ=3$ (the average is $2.20$).
}
\label{fig:dw_2d}
\end{figure}

We also perform a similar analysis for the two-dimensional stochastic differential equation,
\begin{equation}\label{eqn:sde2}
\left\{
\begin{aligned}
dx_t&=y_tdt+\sigma dB_t^x,
\\
dy_t&=[-4y_t+x_t(1-x_t^2)]dt+\sigma dB_t^y,
\end{aligned}
\right.
\end{equation}
where $B_t^x$ and $B_t^y$ denote two independent standard Brownian motions.
Without noise ($\sigma=0$),
there exist two sinks at $(\pm 1,0)$ and one saddle at $(0,0)$.
We set $\sigma=\sqrt{0.08}$ in the following.

We similarly consider two datasets of the form \refeqn{eqn:data_set},
$$
D_1^{\rm 2D}=\{ (\xi^k_n, \eta^k_n)=(x^k_{n\Delta t}, y^k_{n\Delta t}) \in\mathbb{R}^2 \mid k=1,\dots,1000000, \ n=0,1\},
$$
and
$$
D_2^{\rm 2D}=\{ (\xi^k_n, \eta^k_n) \in\mathbb{R}^2 \mid k\in K, \ n=0,1,\dots,99\}
$$
with the index set $K = \{ (a,c) \mid a=1,\dots,30, \ c=0,1,2,3\}$,
by $(\xi^k_n,\eta^k_n)=(x^a_t,y^a_t)$ with $t=0.1n+0.025c$.
Here, for $D_1^{\rm 2D}$,
we take $1000000$ initial points that are
uniformly randomly distributed over the square $[-2,2]\times[-2,2]$, and,
for each initial point, we make one step time evolution
by integrating \refeqn{eqn:sde}
over the time interval $[0,\Delta t]$ with a time step of $\Delta t=0.1$,
whereas for
$D_2^{\rm 2D}$,
we take $30$ initial points uniformly randomly distributed over the square $[-2,2]\times[-2,2]$,
and make 399 steps time evolution by integrating \refeqn{eqn:sde2} over the
time interval $[0, 10-\Delta t]$
with a time step of
$\Delta t=0.025$.

When the multi-valued map is determined from the dataset $D_1^{\rm 2D}$
without using the parameter setting given in \S 2.2, the resulting Morse graph becomes as shown in \reffig{fig:dw_2d}(a),
while it is determined by applying the MGSTD method with the above criterion of parameter setting,
the Morse graph changes to that shown in \reffig{fig:dw_2d}(b),
with two relatively large Morse sets at around $(1,0)$ (orange) and $(-1,0)$ (red),
corresponding to the sink equilibria of the deterministic ODE (without noise),
and a Morse set at around $(0,0)$,
corresponding to the saddle equilibrium.
The result for the dataset $D_2^{\rm 2D}$
is given in \reffig{fig:dw_2d}(c).
In this case, similarly to the one-dimensional case,
the result of Morse graphs depend on the choice of $(\delta_1, \delta_2)$, 
but the MGSTD vector field given in \reffig{fig:dw_2d}(c) qualitatively reproduces the flow 
of \refeqn{eqn:sde2} without noise.

\subsection{Meteorological data}
We here explain the meteorological data to which the MGSTD method is applied.
We obtained a set of time-sequence vectors based on 
a re-analysis dataset, in which 
many types of satellite observation data and special sounding observation are assimilated 
with a weather forecast model and a three-dimensional variational analysis method, 
named JRA25/JCDAS archived by the Japan Meteorological Agency \cite{reanalysis}.
The dataset represents synoptic to global meteorological phenomena on a scale greater 
than several hundred kilometers
with 1.25 degrees by 1.25 degrees mesh spacing, 
and covers a recent period with a 6-hour interval since many meteorological satellites were 
launched.  The analysis period is restricted to three winter months, December, January, and February, from 1979/80 to
2010/11,
and 
then we used
32 (year)
data segments of
$90\mbox{ (day)}\times 4\mbox{ (data per day)}$
length in time. The analysis domain is 
the whole domain north from 20$^\circ$N, with the grid-points being 288 in longitude 
and 57 in latitude. 
After subtracting the trivial seasonal cycle from the data, the low-pass filter extracting 
variations with a period longer than 10 days was taken for the geopotential height anomaly
at a specific isobaric surface of 500 hPa for the tropospheric case and 10 hPa for the stratospheric case.
The isobaric geopotential height is conventionally used by meteorologists for the identification of upper-air low or high pressure systems. 
The principal component analysis applied to the low-frequency variability (LFV) data
eventually provided a set of time-series vectors
with its component being $288\times 57$: 
only the first and second modes that we used explain approximately 25 \% of the LFV 
variance for the tropospheric case and 65 \% of the LFV variance for the stratospheric case, respectively.
The phase space is then spanned by two orthonormal bases of these first and 
second PC modes, just as 
\cite{Inatsu_Nakano_Mukougawa_2013,Inatsu_Nakano_Kusuoka_Mukougawa_2015}.

The dataset thus obtained is,
for both the troposphere and the stratosphere,
$$
\Tilde{D}^{\rm L} = \{ (\tilde{y}_1(a,b,c),\tilde{y}_2(a,b,c)) \in\mathbb{R}^2 \mid
a=1,\dots, 32, \ b=1,\dots, 90, \ c=0,6,12,18 \}
$$
where ${\rm L}$ indicates either the troposphere or the stratosphere,
$\tilde{y}_\ell(a,b,c)$ $(\ell=1,2)$ is the score of the PC$\ell$,
and $a$, $b$, and $c$ denote year, day, and o'clock, respectively.
We consider this as the dataset of the form \refeqn{eqn:data_set},
$$
D^{\rm L} = \{ (y_1^k(n),y_2^k(n)) \in\mathbb{R}^2 \mid k\in K, \ n=1,\dots, 90 \}
$$
with the index set $K = \{ (a,c) \mid a=1,\dots,32, \ c=0,6,12,18\}$,
by $y_\ell^k(n)=\tilde{y}_\ell(a,n,c)$;
notice that $(y_1^k(n), y_2^k(n))$ corresponds to $y^k_n$ in \refeqn{eqn:data_set}.
Accordingly, the subset of the form \refeqn{eqn:data_subset} is,
$$
D'^{\rm L} = \{ (y_1^k(n),y_2^k(n)) \in\mathbb{R}^2 \mid k\in K, \ n=1,\dots, 89 \}
$$
The MGSTD method is hence applied to this dataset $D^{\rm L}$.

\begin{figure}[tbhp]
\begin{minipage}[t]{0.49\textwidth}
\centering
(a) Troposphere\\
\includegraphics[width=\textwidth]{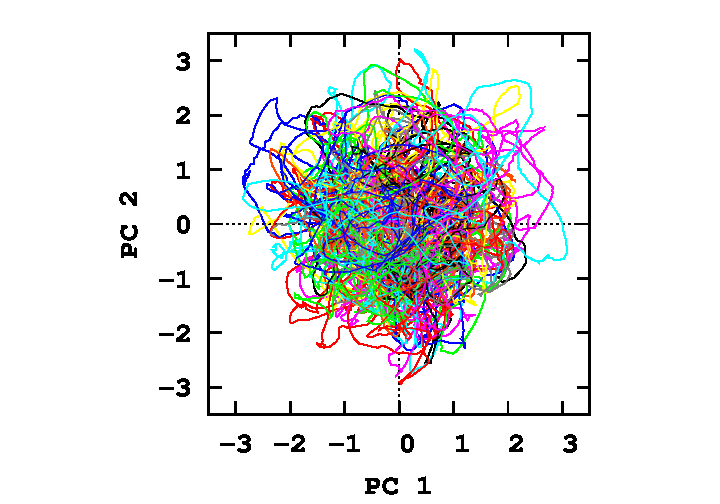}
\end{minipage}
\begin{minipage}[t]{0.49\textwidth}
\centering
(b) Stratosphere\\
\includegraphics[width=\textwidth]{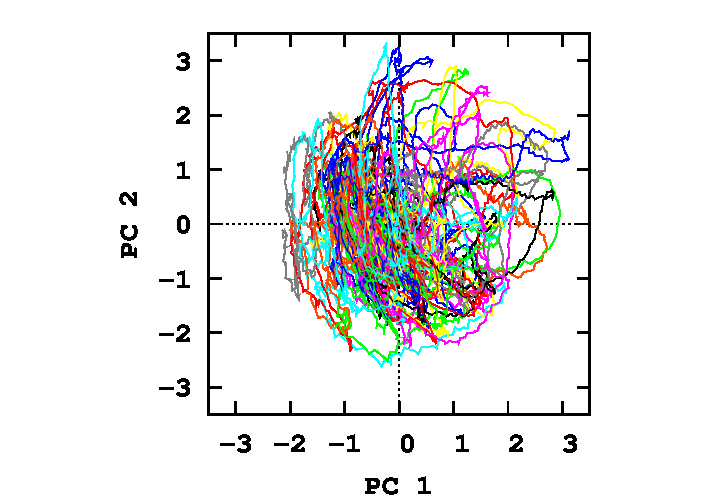}
\end{minipage}
\caption{
Trajectories projected onto the PC1-PC2 plane for (a) the troposphere and (b) the stratosphere.
Different colors denote different years.
}
\label{fig:traj}
\end{figure}

\begin{figure}[tbhp]
\begin{minipage}[t]{0.49\textwidth}
\centering
(a) Troposphere\\
\includegraphics[width=\textwidth]{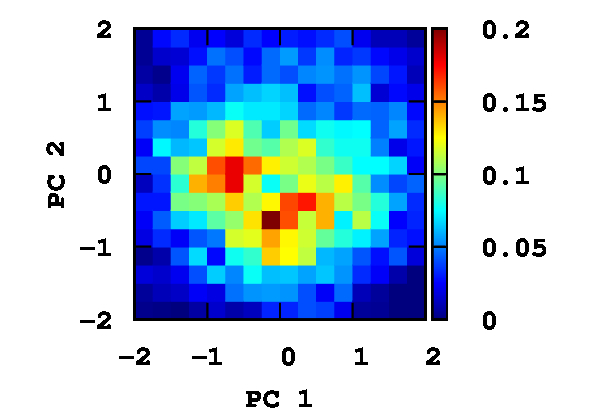}
\end{minipage}
\begin{minipage}[t]{0.49\textwidth}
\centering
(b) Stratosphere\\
\includegraphics[width=\textwidth]{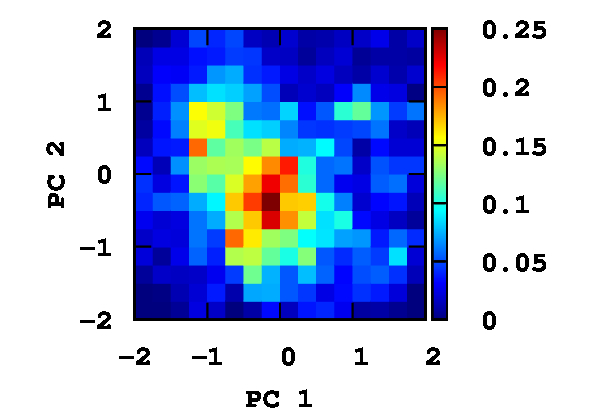}
\end{minipage}
\caption{
Probability density functions of PC1 and PC2 for (a) the troposphere and (b) the stratosphere. The grid size $h=0.25$.
}
\label{fig:pdf2}
\end{figure}

The trajectories of thus obtained time-series data of the PC scores
for the troposphere and the stratosphere
are projected onto the PC1-2 plane in \reffig{fig:traj},
which show highly stochastic dynamics.
The probability density functions of PC1 and PC2,
shown in \reffig{fig:pdf2},
display skewed, non-Gaussian forms.
This is related to the existence of several known, persistent,
characteristic patterns of the pressure field
detected as local departure from two-dimensional Gaussian PDF in the phase space.
For the troposphere (\reffig{fig:pdf2}-(a)), in particular,
those at around (PC1, PC2) 
$= (0, -1)$, $(1,0)$, $(1,1)$, and $(-1,0)$
are called ZNAO, PNA, 
BNAO, and RNA, respectively.
(Note that the PC2 is upside down compared with \cite{Kimoto_Ghil_1993a}.) 
 
These spatial patterns are dominant week-to-month variability and have a great impact to weather systems in the Northern Hemisphere. For example, a pair of ZNAO and BNAO are long recognized as a see-saw pattern of Azores high and Icelandic low pressures in the North Atlantic, which is usually called the North Atlantic Oscillation (NAO)
(originally \cite{Walker_1924} but see \cite{Hurrell_et_al_2003} for the overview of NAO studies). The positive phase of NAO, say ZNAO, shifts the jet stream, the storm-track, and precipitation poleward, while the negative phase of NAO shifts them equatorward. On the other hand, a pair of PNA and RNA has been recognized as a see-saw pattern across the Pacific to North America, which is usually called the Pacific-North American (PNA) pattern. A positive phase with intensifying the Aleutian low and a negative phase respectively corresponds to PNA and RNA in this paper. The typical Rossby-wave train propagates from the equatorial North Pacific to the south-eastern US via Alaska \cite{Wallace_Gutzler_1981}. A transition between ZNAO and RNA and another transition cycling PNA, BNAO, RNA, ZNAO, and backing to PNA were statistically discussed 
in \cite{Kimoto_Ghil_1993b}.

\section{Result}

\begin{figure}[tbhp]
\begin{minipage}[t]{0.49\textwidth}
\centering
(a) Troposphere\\
\includegraphics[width=\textwidth]{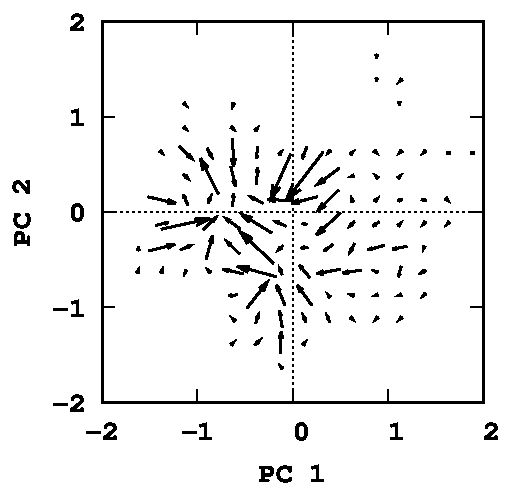}
\end{minipage}
\begin{minipage}[t]{0.49\textwidth}
\centering
(b) Stratosphere\\
\includegraphics[width=\textwidth]{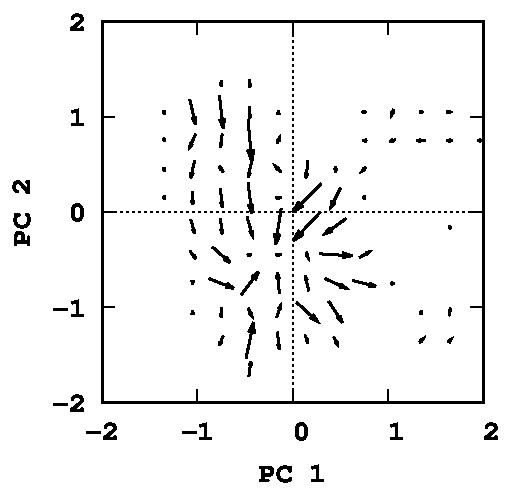}
\end{minipage}
\caption{
MGSTD vector fields for
(a) the troposphere with the dataset $D^{\rm tropo}$, and
$h=0.25$, $\rho=1.1$, $\mu_*^\circ=8$ (the average $7.76$),
and
(b) the stratosphere with the dataset $D^{\rm strato}$, and
$h=0.3$, $\rho=1.1$, $\mu_*^\circ=15$ (the average $14.51$).
}
\label{fig:result}
\end{figure}

\begin{figure}[htbp]
\centering
\tabcolsep=0mm
\begin{tabular}{ccc}
\begin{minipage}[c]{0.3\textwidth}
\includegraphics[width=\textwidth]{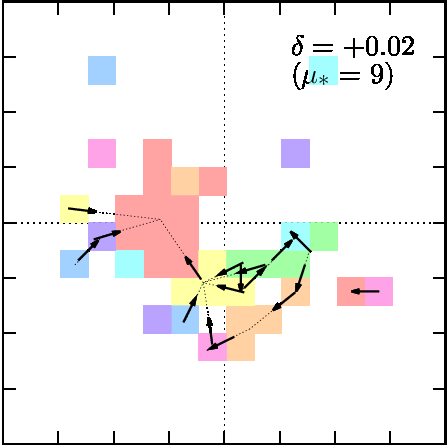}
\end{minipage}
&
\begin{minipage}[c]{0.3\textwidth}
\includegraphics[width=\textwidth]{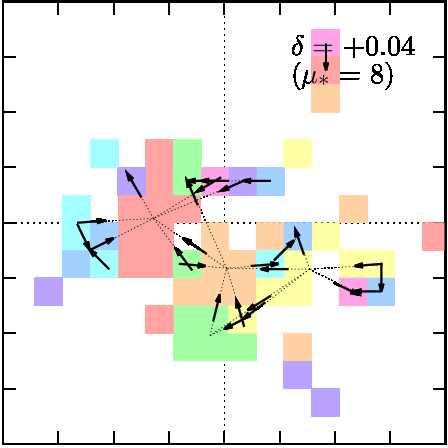}
\end{minipage}
&
\begin{minipage}[c]{0.3\textwidth}
\includegraphics[width=\textwidth]{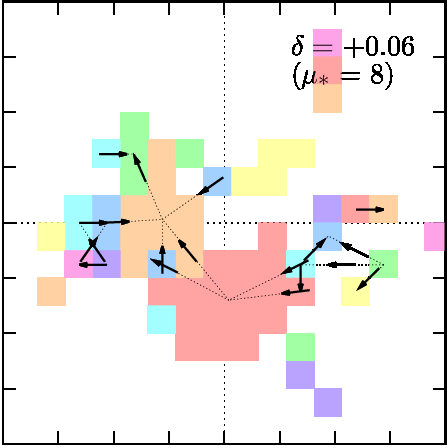}
\end{minipage}
\\
\begin{minipage}[c]{0.3\textwidth}
\includegraphics[width=\textwidth]{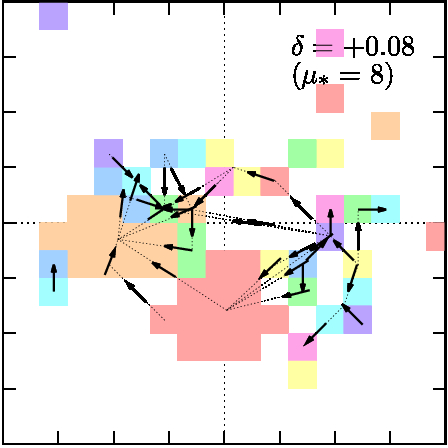}
\end{minipage}
&
\begin{minipage}[c]{0.3\textwidth}
\includegraphics[width=\textwidth]{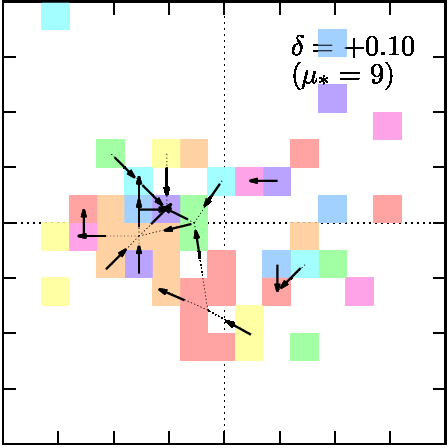}
\end{minipage}
&
\begin{minipage}[c]{0.3\textwidth}
\includegraphics[width=\textwidth]{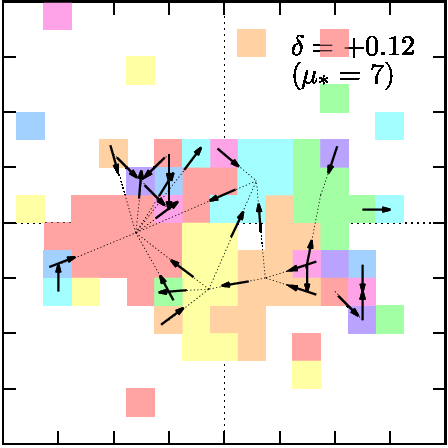}
\end{minipage}
\\
\begin{minipage}[c]{0.3\textwidth}
\includegraphics[width=\textwidth]{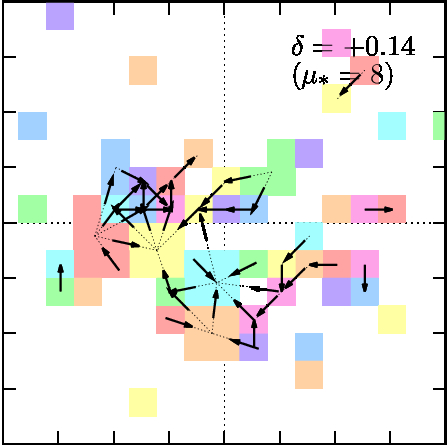}
\end{minipage}
&
\begin{minipage}[c]{0.3\textwidth}
\includegraphics[width=\textwidth]{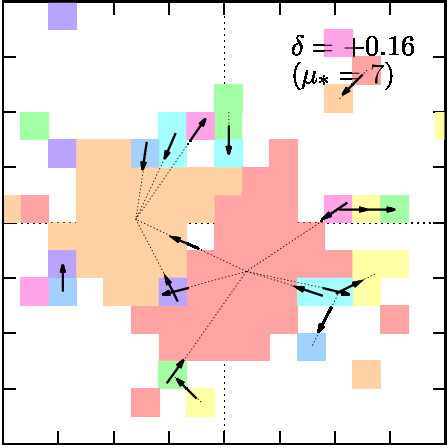}
\end{minipage}
&
\begin{minipage}[c]{0.3\textwidth}
\includegraphics[width=\textwidth]{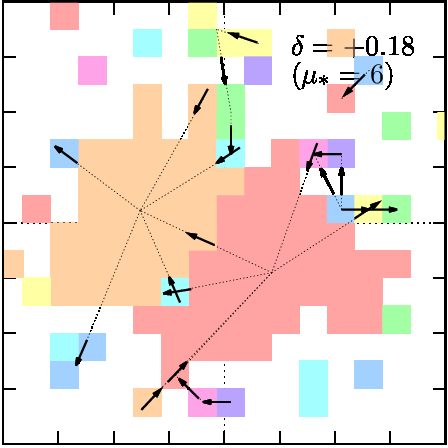}
\end{minipage}
\\
\begin{minipage}[c]{0.3\textwidth}
\includegraphics[width=\textwidth]{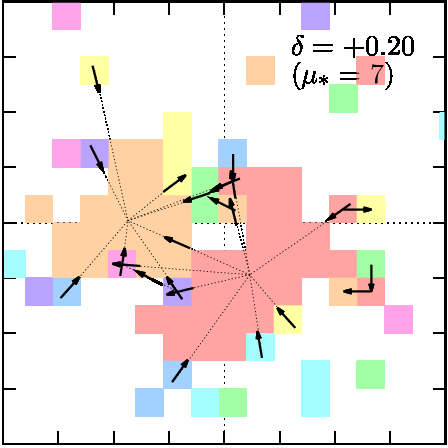}
\end{minipage}
&
\begin{minipage}[c]{0.3\textwidth}
\includegraphics[width=\textwidth]{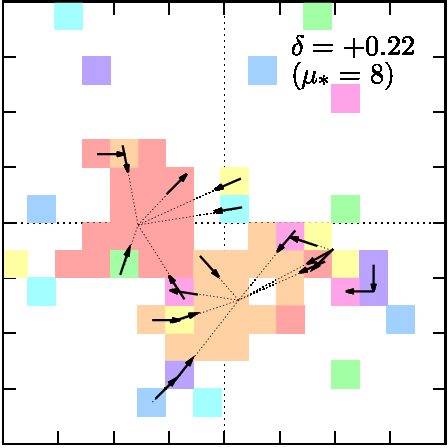}
\end{minipage}
&
\begin{minipage}[c]{0.3\textwidth}
\includegraphics[width=\textwidth]{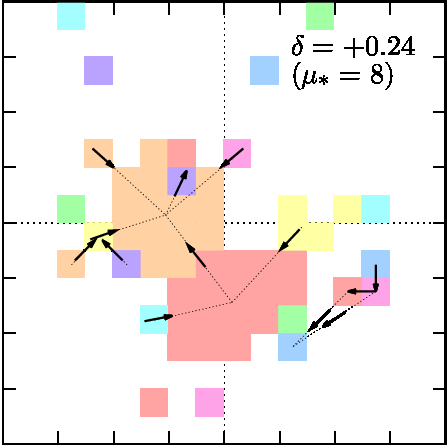}
\end{minipage}
\end{tabular}
\caption{
Morse sets and the phase space presentation of the Morse graphs for $D^{\rm tropo}$ 
with $(\delta_1, \delta_2)=\delta(1,0)$ for varying
$\delta=0.02$ to $0.24$
(see \reffig{fig:Morse_graph_set_tropo} for $\delta=0$);
$h=0.25$, $\rho=1.1$.
The abscissa and ordinate are PC1 and PC2, respectively, in all panels.
}
\label{fig:MS_vs_delta_tropo}
\end{figure}

To determine the parameters for use in the analysis,
as discussed in \S 2.2, we first chose $m=2$, and studied the datasets $D^{\rm tropo}$ for the troposphere and 
$D^{\rm strato}$ for the stratosphere.  We chose $h=0.25$ from $D^{\rm tropo}$ and
$h=0.3$ from $D^{\rm strato}$,
according to the criterion from \reffig{fig:hselect} and
in the accompanying explanation.

Fixing $\rho=1.1$,
we computed $\mu_*^{\circ}$ by varying $\delta_1, \delta_2$.
Its average values over $\delta_1,\delta_2$ are $\mu_*^{\circ}=7.76$ for the troposphere
and $\mu_*^{\circ}=14.51$ for the stratosphere.  Since $\mu_*$ must be an integer by definition,
we chose $\mu_*^{\circ}=8$ for the troposphere and $\mu_*^{\circ}=15$
for the stratosphere.

\medskip

For the troposphere, the MGSTD vector field for $D^{\rm tropo}$ with $h=0.25$, 
$\rho=1.1$, and $\mu_*^\circ=8$
is shown in \reffig{fig:result}(a).
We observe leftward and down-leftward motions in the fourth quadrant,
and up-leftward motion in the third quadrant.
We also observe a motion from the first to the second quadrant through a vicinity of the origin.

The phase space presentations of the Morse graphs for various choices of 
$\delta_1,\delta_2$ are shown in \reffig{fig:MS_vs_delta_tropo},
together with the corresponding $\mu_*^{\circ}$.
In most cases,
the motions observed in the MGSTD vector field
are consistent with transitions between Morse sets.
In particular,
the up-leftward motion in the third quadrant
seems to correspond to the transition between the two largest Morse sets,
from the one at around $(0,-1)$ to that at around $(-1,0)$.
Moreover, in some cases,
a Morse set around $(1,0)$ seems to be involved in
the leftward and down-leftward motions in the fourth quadrant,
especially as a (successive) transition to the Morse set at around $(0,-1)$,
and in the motion from the first to the second quadrant through a vicinity of the origin.

\begin{figure}[htbp]
\tabcolsep=0mm
\begin{tabular}{c|ccc}
\backslashbox{$h$}{$\rho$} & 1.1 & 1.3 & 1.5
\\ \hline
\begin{minipage}[c]{0.07\textwidth}
$0.20$
\end{minipage}
&
\begin{minipage}[c]{0.31\textwidth}
\includegraphics[width=\textwidth]{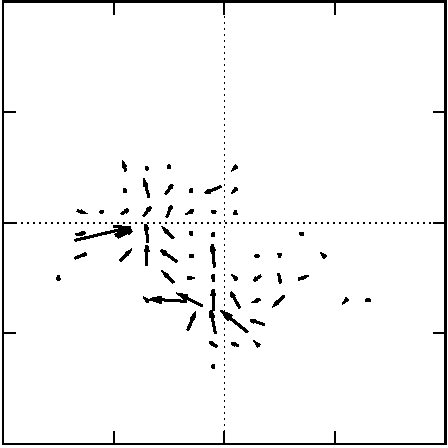}
\end{minipage}
&
\begin{minipage}[c]{0.31\textwidth}
\includegraphics[width=\textwidth]{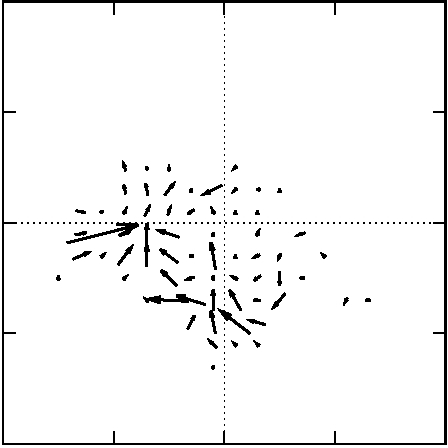}
\end{minipage}
&
\begin{minipage}[c]{0.31\textwidth}
\includegraphics[width=\textwidth]{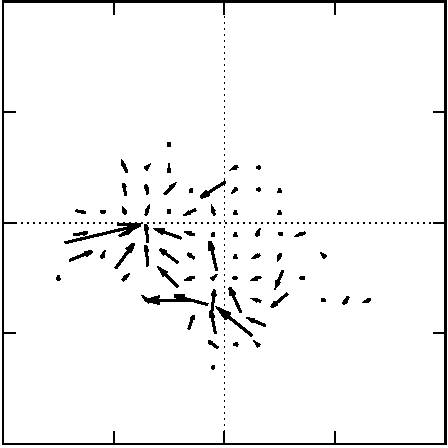}
\end{minipage}
\\
\begin{minipage}[c]{0.07\textwidth}
$0.25$
\end{minipage}
&
\begin{minipage}[c]{0.31\textwidth}
\includegraphics[width=\textwidth]{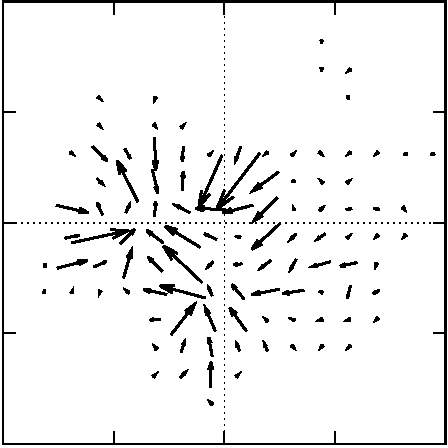}
\end{minipage}
&
\begin{minipage}[c]{0.31\textwidth}
\includegraphics[width=\textwidth]{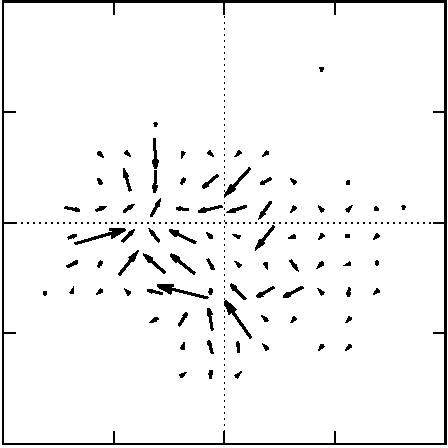}
\end{minipage}
&
\begin{minipage}[c]{0.31\textwidth}
\includegraphics[width=\textwidth]{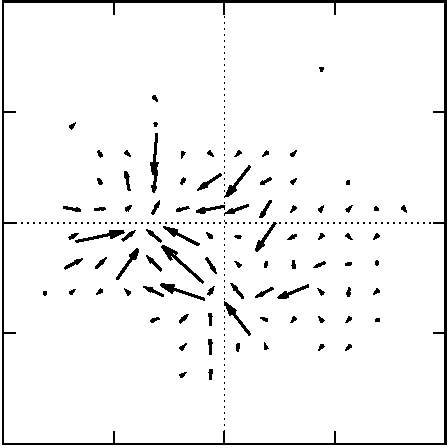}
\end{minipage}
\\
\begin{minipage}[c]{0.07\textwidth}
$0.30$
\end{minipage}
&
\begin{minipage}[c]{0.31\textwidth}
\includegraphics[width=\textwidth]{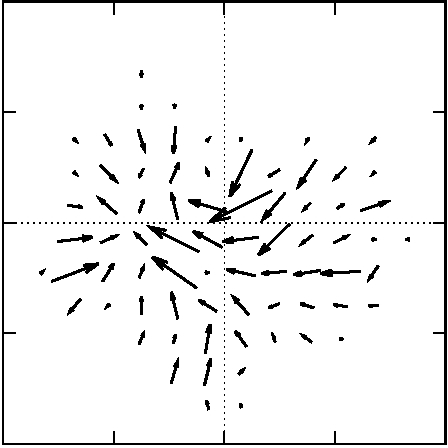}
\end{minipage}
&
\begin{minipage}[c]{0.31\textwidth}
\includegraphics[width=\textwidth]{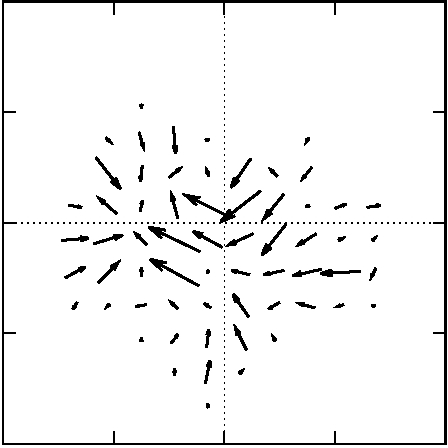}
\end{minipage}
&
\begin{minipage}[c]{0.31\textwidth}
\includegraphics[width=\textwidth]{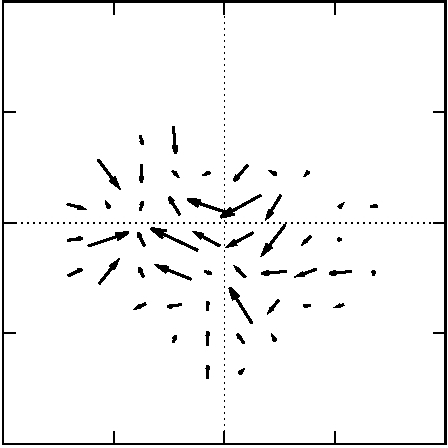}
\end{minipage}
\end{tabular}
\caption{
MGSTD vector fields for the troposphere with $D^{\rm tropo}$
for varying $h$ and $\rho$.
The values of $\mu_*^\circ$ (resp. the average) are, from left to right, 
top to bottom, 
6 (5.29), 6 (5.36), 6 (5.42), 
8 (7.76), 9 (8.05), 9(8.23), 
11 (10.52), 12 (11.50), 13 (12.17).  
The abscissa and ordinate are PC1 and PC2, respectively, in all panels.
}
\label{fig:morse_arrow_avg_tropo_1}
\end{figure}

These dominant motions are commonly observed over a finite range of $h$ and $\rho$,
as shown in \reffig{fig:morse_arrow_avg_tropo_1}.

\medskip

For the stratosphere, the MGSTD vector field for $D^{\rm strato}$ with $h=0.3$, $\rho=1.1$,
and $\mu_*^\circ=15$
is shown in \reffig{fig:result}(b).
We observe motions to a sink-like location at around $(-0.25, -0.25)$,
and a downward motion in the second and third quadrants.
We also observe a leftward motion in the first quadrant,
though it is less prominent than the previous one.

\begin{figure}[htbp]
\centering
\tabcolsep=0mm
\begin{tabular}{ccc}
\begin{minipage}[c]{0.3\textwidth}
\includegraphics[width=\textwidth]{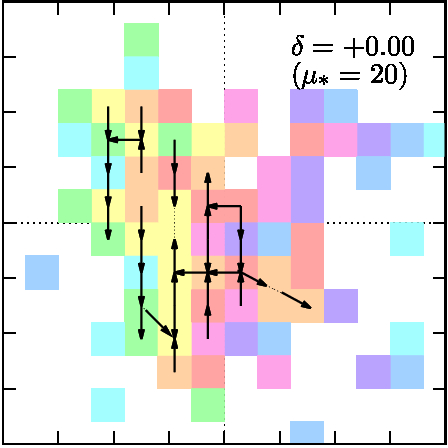}
\end{minipage}
&
\begin{minipage}[c]{0.3\textwidth}
\includegraphics[width=\textwidth]{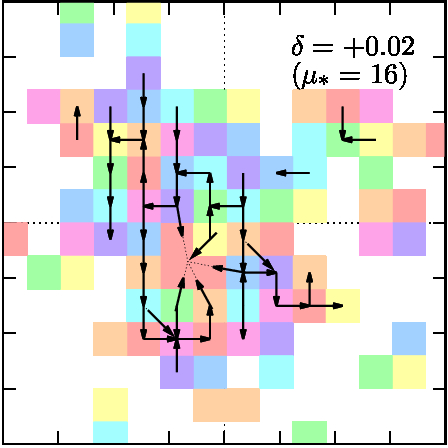}
\end{minipage}
&
\begin{minipage}[c]{0.3\textwidth}
\includegraphics[width=\textwidth]{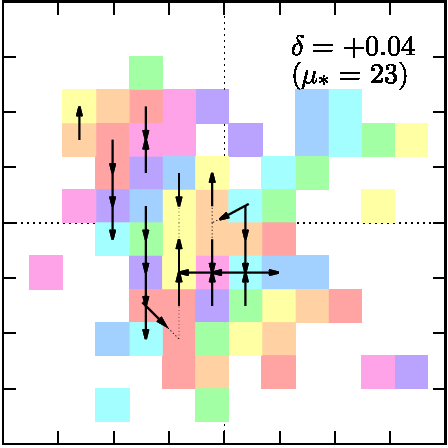}
\end{minipage}
\\
\begin{minipage}[c]{0.3\textwidth}
\includegraphics[width=\textwidth]{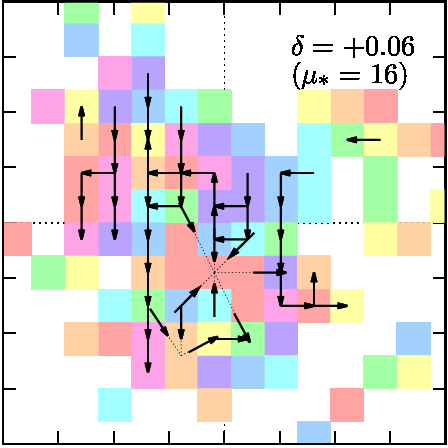}
\end{minipage}
&
\begin{minipage}[c]{0.3\textwidth}
\includegraphics[width=\textwidth]{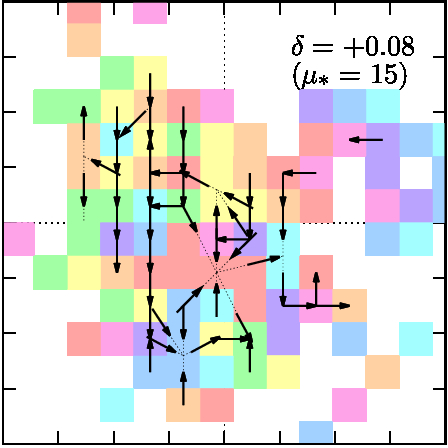}
\end{minipage}
&
\begin{minipage}[c]{0.3\textwidth}
\includegraphics[width=\textwidth]{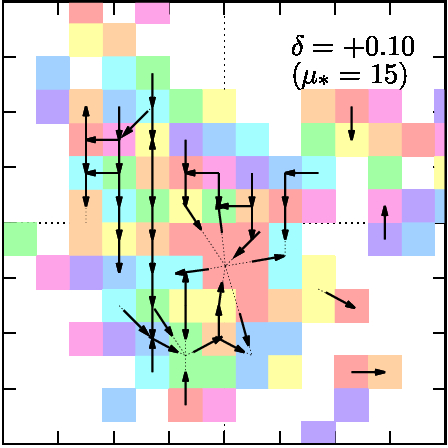}
\end{minipage}
\\
\begin{minipage}[c]{0.3\textwidth}
\includegraphics[width=\textwidth]{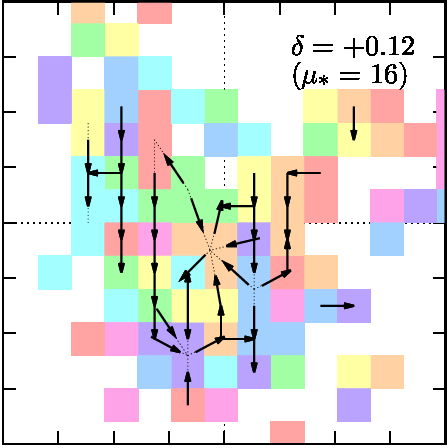}
\end{minipage}
&
\begin{minipage}[c]{0.3\textwidth}
\includegraphics[width=\textwidth]{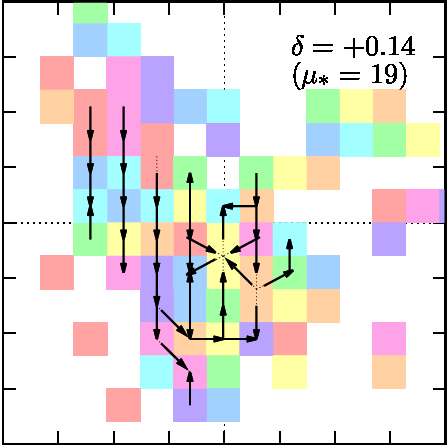}
\end{minipage}
&
\begin{minipage}[c]{0.3\textwidth}
\includegraphics[width=\textwidth]{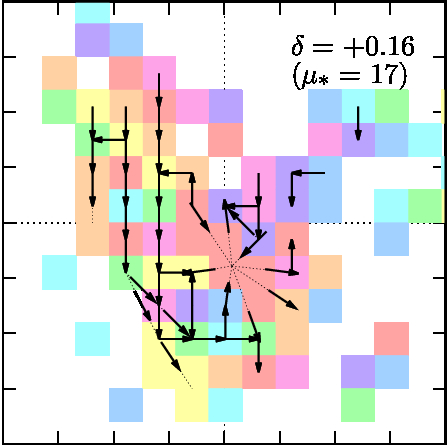}
\end{minipage}
\\
\begin{minipage}[c]{0.3\textwidth}
\includegraphics[width=\textwidth]{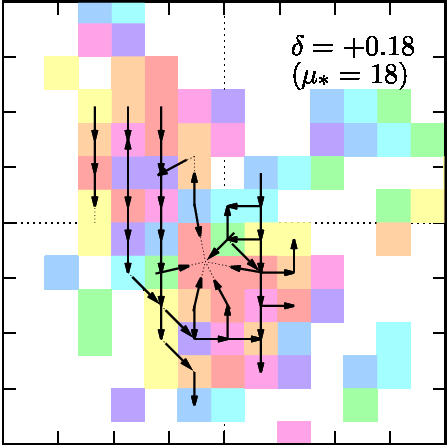}
\end{minipage}
&
\begin{minipage}[c]{0.3\textwidth}
\includegraphics[width=\textwidth]{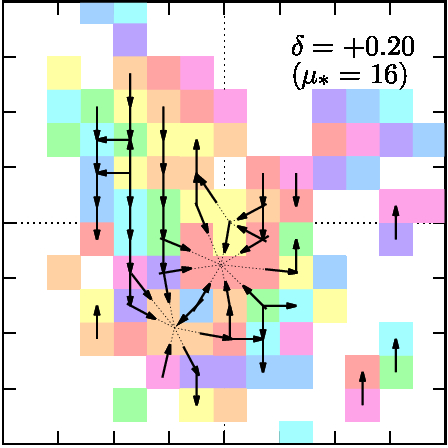}
\end{minipage}
&
\begin{minipage}[c]{0.3\textwidth}
\includegraphics[width=\textwidth]{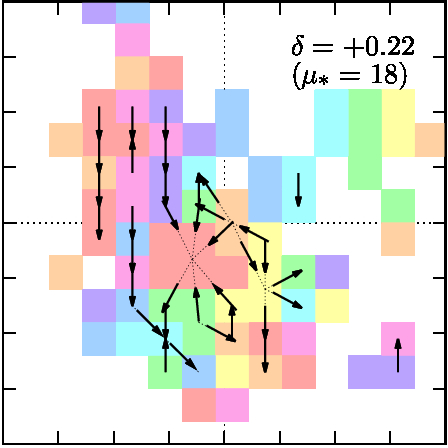}
\end{minipage}
\\
\begin{minipage}[c]{0.3\textwidth}
\includegraphics[width=\textwidth]{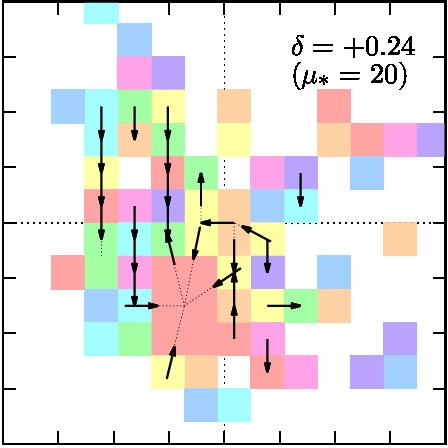}
\end{minipage}
&
\begin{minipage}[c]{0.3\textwidth}
\includegraphics[width=\textwidth]{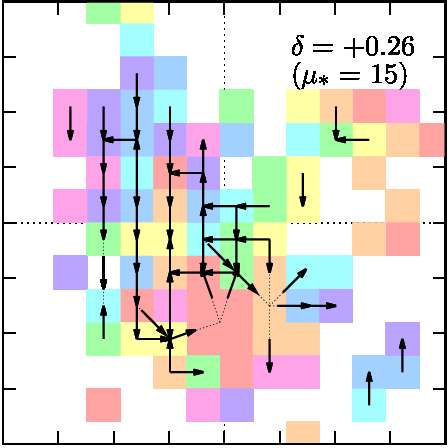}
\end{minipage}
&
\begin{minipage}[c]{0.3\textwidth}
\includegraphics[width=\textwidth]{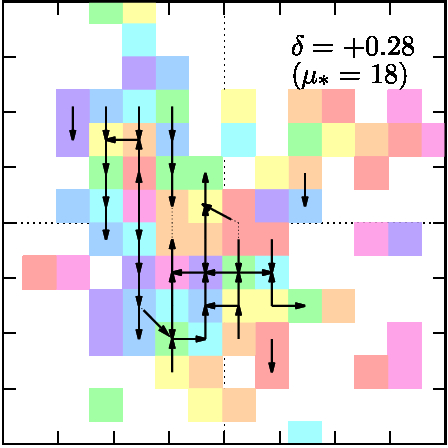}
\end{minipage}
\end{tabular}
\caption{
Morse sets and the phase space presentation of the Morse graphs for $D^{\rm strato}$ 
with $(\delta_1, \delta_2)=\delta(1,0)$ for varying $\delta=0.0$ to $0.28$; 
$h=0.3$, $\rho=1.1$.
The abscissa and ordinate are PC1 and PC2, respectively, in all panels.
}
\label{fig:MS_vs_delta_strato}
\end{figure}

The phase space presentations of the Morse graphs for various choices of 
$\delta_1,\delta_2$ are shown in \reffig{fig:MS_vs_delta_strato},
together with the corresponding $\mu_*^{\circ}$.
In most cases, the motions observed in the MGSTD vector field are again
consistent with transitions between Morse sets.
In particular, the sink seems to 
correspond to a relatively large Morse set at around $(-0.25, -0.25)$
that tends to attract transitions inward.
Moreover, the downward motion in the second and third quadrants
seems to correspond to successive transitions between Morse sets downward in the second and third quadrants.
In addition, in some cases, the leftward motion in the first quadrant
seems to correspond to the leftward or left-and-downward transitions between Morse sets in the first quadrant.

\begin{figure}[htbp]
\tabcolsep=0mm
\begin{tabular}{c|ccc}
\backslashbox{$h$}{$\rho$} & 1.1 & 1.3 & 1.5
\\ \hline
\begin{minipage}[c]{0.07\textwidth}
$0.25$
\end{minipage}
&
\begin{minipage}[c]{0.31\textwidth}
\includegraphics[width=\textwidth]{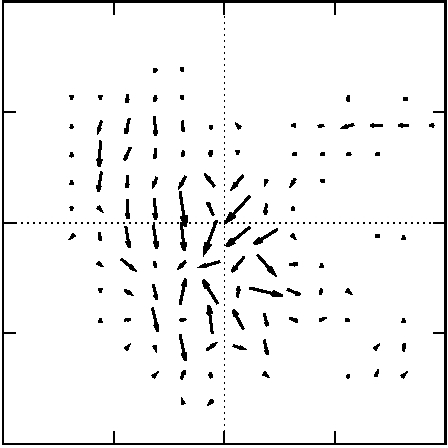}
\end{minipage}
&
\begin{minipage}[c]{0.31\textwidth}
\includegraphics[width=\textwidth]{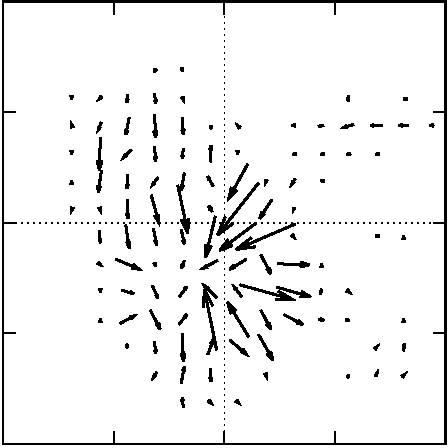}
\end{minipage}
&
\begin{minipage}[c]{0.31\textwidth}
\includegraphics[width=\textwidth]{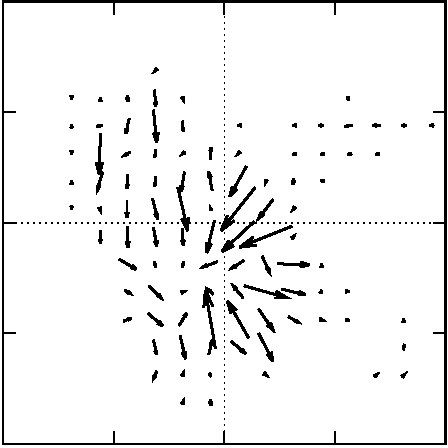}
\end{minipage}
\\
\begin{minipage}[c]{0.07\textwidth}
$0.30$
\end{minipage}
&
\begin{minipage}[c]{0.31\textwidth}
\includegraphics[width=\textwidth]{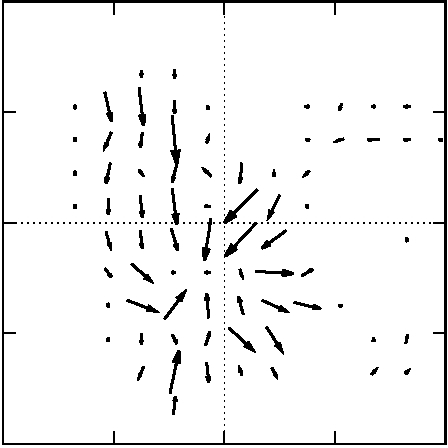}
\end{minipage}
&
\begin{minipage}[c]{0.31\textwidth}
\includegraphics[width=\textwidth]{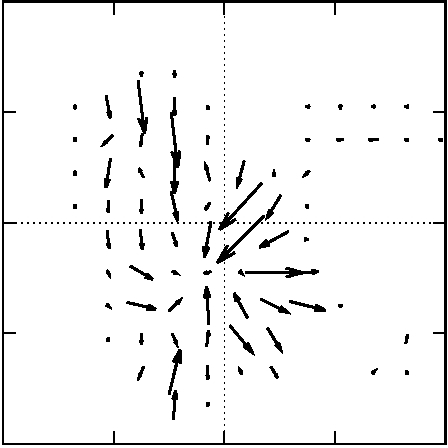}
\end{minipage}
&
\begin{minipage}[c]{0.31\textwidth}
\includegraphics[width=\textwidth]{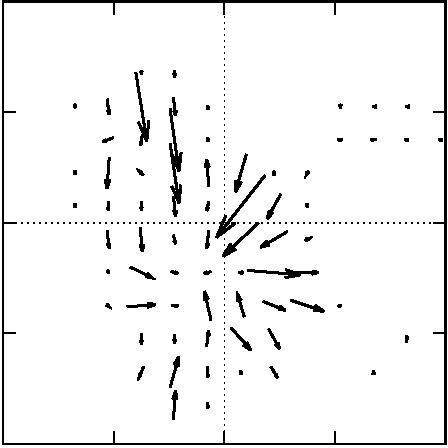}
\end{minipage}
\\
\begin{minipage}[c]{0.07\textwidth}
$0.35$
\end{minipage}
&
\begin{minipage}[c]{0.31\textwidth}
\includegraphics[width=\textwidth]{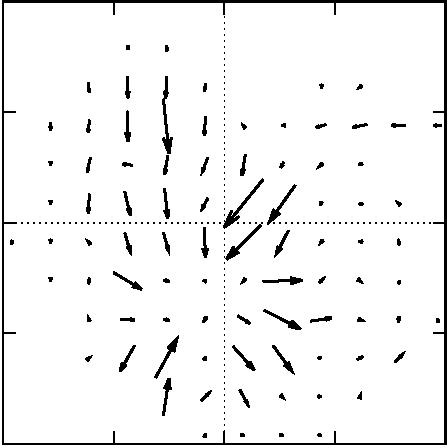}
\end{minipage}
&
\begin{minipage}[c]{0.31\textwidth}
\includegraphics[width=\textwidth]{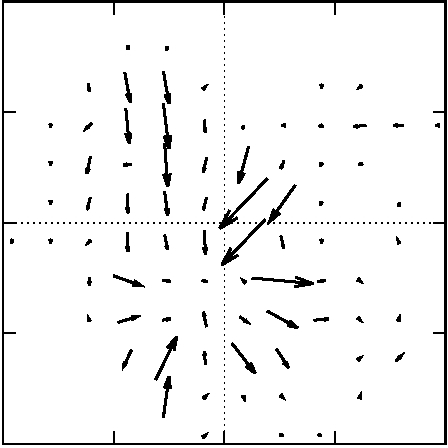}
\end{minipage}
&
\begin{minipage}[c]{0.31\textwidth}
\includegraphics[width=\textwidth]{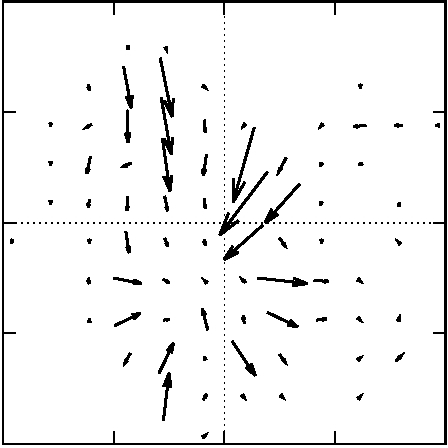}
\end{minipage}
\end{tabular}
\caption{
MGSTD vector field for the stratosphere with $D^{\rm strato}$
for varying $h$ and $\rho$.
The values of $\mu_*^\circ$ (resp. the average) are, from left to right, 
top to bottom, 
13 (12.06), 13 (12.94), 14 (13.60),
15 (14.51), 16 (15.93), 17 (16.51),
16 (15.90), 18 (17.80), 19 (18.82).
The abscissa and ordinate are PC1 and PC2, respectively, in all panels.
}
\label{fig:morse_arrow_avg_strato_1}
\end{figure}

These dominant motions are again commonly observed over a finite range of $h$ and $\rho$,
as shown in \reffig{fig:morse_arrow_avg_strato_1}.

\section{Discussion and concluding remarks}

For the troposphere, we have observed two dominant motions.
One is left-downward in the fourth quadrant,
and the other is left-upward in the third quadrant,
by detecting
the clockwise circular dynamics from a Morse set around $(1,0)$ to one around $(0,-1)$ 
and then to another around $(-1,0)$. 
Noting that the choice of signs of PC scores is arbitrary,
a comparison of the probability density function in \reffig{fig:pdf2}(a) with that reported previously in [16]
indicates that these motions correspond to the transitions
from PNA to ZNAO and that from ZNAO to RNA, respectively.
Luo et al. \cite{Luo_Cha_Feldstein} investigated
a characteristic weather pattern in a transition path from ZNAO to BNAO
and another transition path from BNAO and ZNAO in the Atlantic Ocean.
The former transition is through the Scandinavia blocking high with a negative PC1 projection,
and the latter is through the Atlantic ridge with a positive PC1 projection \cite{Woolings_Pinto_Santos}.
This is consistent with our analysis. 

For the stratosphere,
we have observed a dominant motion downward from the second to the third quadrant.
This is associated with a wave-energy charge by the vertical propagation of Rossby waves with their zonal wavenumber one.
We have also observed another dominant sink-like location at around $(-0.25, -0.25)$.
This corresponds to the resetting of the characteristic pattern,
as PC1 and PC2 both become near zero.
The small difference results from the skewness of the data set distribution.
There is also a less prominent motion directed
leftward in the first quadrant.
This may be related to the final stage of an event called the stratospheric sudden warming.
These results are also consistent with earlier studies
\cite{Inatsu_Nakano_Kusuoka_Mukougawa_2015}.

For the troposphere,
the transitions between weather regimes which previous studies detected
have relied on meteorological knowledge.
An advantage of our proposed method is that such patterns can be identified without any meteorological foresight.
For the stratosphere, the stratospheric sudden warming
and the subsequent polar-vortex amplification
were well known and explained with atmospheric dynamics.
However, the method presented in this paper can identify this transition only from the time-series data.

The results demonstrate the efficacy of the proposed method
not only for time-series derived from model dynamical systems
but also for that from real measurements.
Furthermore, this kind of automatic analysis techniques,
including machine learning,
could be helpful in scientific fields with a massive amount of data like meteorology.

\bigskip

\noindent{\bf Acknowledgement}\\
The authors are grateful to Masaki Nomura for helping to calculate the Morse graphs.
This work was supported by JST CREST, and by JSPS KAKENHI Grant Numbers JP25287029,  JP26310208, and JP18H03671.
MI was partly supported by JSPS KAKENHI Grant Numbers JP25610028, JP26310201, and JP18K03734.

\bibliographystyle{siam}
\bibliography{script_v8.bbl}

\begin{thebibliography}{10}

\bibitem{AUTO}
{\em Auto --- software for continuation and bifurcation problems in ordinary
  differential equations}.
\newblock \texttt{http://indy.cs.concordia.ca/auto/}.

\bibitem{CAPD}
{\em Computer assisted proofs in dynamics}.
\newblock \texttt{http://capd.wsb-nlu.edu.pl/}.

\bibitem{Aeyels_1981}
{\sc Dirk Aeyels}, {\em Generic observability of differentiable systems}, SIAM
  Journal on Control and Optimization, 19 (1981), pp.~595--603.

\bibitem{Arai_etal_2009}
{\sc Z.~{Arai}, W.~{Kalies}, H.~{Kokubu}, K.~{Mischaikow}, H.~{Oka}, and
  P.~{Pilarczyk}}, {\em {A Database Schema for the Analysis of Global Dynamics
  of Multiparameter Systems}}, SIAM Journal on Applied Dynamical Systems, 8
  (2009), pp.~757--789.

\bibitem{Ban_Kalies_2006}
{\sc Hyunju Ban and William~D Kalies}, {\em A computational approach to
  Conley's decomposition theorem}, Journal of Computational and Nonlinear
  Dynamics, 1 (2006), pp.~312--319.


\bibitem{Bush_etal_2012}
{\sc J.~{Bush}, M.~{Gameiro}, S.~{Harker}, H.~{Kokubu}, K.~{Mischaikow},
  I.~{Obayashi}, and P.~{Pilarczyk}}, {\em {Combinatorial-topological framework
  for the analysis of global dynamics}}, Chaos, 22 (2012), p.~047508.

\bibitem{Conley_1978}
{\sc Charles~C. Conley}, {\em Isolated invariant sets and the Morse index},
  no.~38 in CBMS Regional Conference Series in Mathematics, American
  Mathematical Society, Providence, R. I., 1978.

\bibitem{GAIO}
{\sc Michael Dellnitz, Gary Froyland, and Oliver Junge}, {\em The algorithms
  behind gaio -- set oriented numerical methods for dynamical systems}, in
  Ergodic theory, analysis, and efficient simulation of dynamical systems,
  Bernold Fiedler, ed., Springer, 2001, pp.~145--174.

\bibitem{Holton_Hakim_2012}
{\sc James~R Holton and Gregory~J Hakim}, {\em An introduction to dynamic
  meteorology}, Academic Press, 2012.
  
\bibitem{Honeycutt1992}
{\sc Rebecca~L~Honeycurtt}, {\em Stochastic Runge-Kutta algorithms, I: White noise},
Physical Review A, 45 (1992), pp.~600-603.

\bibitem{Hurrell_et_al_2003}
{\sc J.~W. Hurrell, Y.~Kushnir, G.~Ottersen, and M.~Visbeck}, {\em {An overview 
of the North Atlantic Oscillation. The North Atlantic Oscillation: Climatic 
Significance and Environmental Impact}}, Geophys. Monogr., 134 (2003), Amer. 
Geophys. Union, pp.~1-36.

\bibitem{Inatsu_Nakano_Kusuoka_Mukougawa_2015}
{\sc M.~{Inatsu}, N.~{Nakano}, S.~{Kusuoka}, and H.~{Mukougawa}}, {\em
  {Predictability of wintertime stratospheric circulation examined by
  non-stationary fluctuation dissipation relation}}, Journal of Atmospheric
  Sciences, 72 (2015), pp.~774--786.

\bibitem{Inatsu_Nakano_Mukougawa_2013}
{\sc M.~{Inatsu}, N.~{Nakano}, and H.~{Mukougawa}}, {\em {Dynamics and
  practical predictability of extratropical wintertime low-frequency
  variability in a low-dimensional system}}, Journal of Atmospheric Sciences,
  70 (2013), pp.~939--952.

\bibitem{Kalies_Mischaikow_Vandervorst_2005}
{\sc William~D Kalies, Konstantin Mischaikow, and Robert~CAM Vandervorst}, {\em
  An algorithmic approach to chain recurrence}, Foundations of Computational
  Mathematics, 5 (2005), pp.~409--449.

\bibitem{Kalney_2002}
{\sc Eugenia Kalney}, {\em Atmospheric modeling, data assimilation and
  predictability}, Cambridge University Press, 2002.

\bibitem{Kimoto_Ghil_1993a}
{\sc M.~{Kimoto} and M.~{Ghil}}, {\em {Multiple Flow Regimes in the Northern
  Hemisphere Winter. Part I: Methodology and Hemispheric Regimes.}}, Journal of
  Atmospheric Sciences, 50 (1993), pp.~2625--2644.

\bibitem{Kimoto_Ghil_1993b}
\leavevmode\vrule height 2pt depth -1.6pt width 23pt, {\em {Multiple Flow
  Regimes in the Northern Hemisphere Winter. Part II: Sectorial Regimes and
  Preferred Transitions.}}, Journal of Atmospheric Sciences, 50 (1993),
  pp.~2645--2673.

\bibitem{Kokubu_Morita_Nomura_Obayashi}
{\sc Hiroshi {Kokubu}, Hidetoshi {Morita}, Masaki {Nomura}, and Ippei
  {Obayashi}}, {\em {Conley-Morse graph analysis of time series}}.
\newblock in preparation.

\bibitem{Luo_Cha_Feldstein}
{\sc D. Luo, J. Cha, and S.B. Feldstein}, {\em Weather Regime Transitions and the 
Interannual Variability of the North Atlantic Oscillation. Part II: Dynamical Processes.} 
J. Atmos. Sci., 69 (2012), pp. 2347--2363.

\bibitem{Miyaji_etal_2016}
{\sc T. {Miyaji}, P. {Pilarczyk}, M. {Gameiro}, H. {Kokubu}, K. {Mischaikow}},
{\em A study of rigorous ODE integrators for multi-scale set-oriented computations},
Applied Numerical Mathematics, 107 (2016), pp.~34--47.

\bibitem{Oksendal}
{\sc Bernt {\O}ksendal}, {\em Stochastic differential equations}, Springer,
  sixth~ed., 2013.

\bibitem{reanalysis}
{\sc Kazutoshi Onogi, Junichi Tsutsui, Hiroshi Koide, Masami Sakamoto, Shinya
  Kobayashi, Hiroaki Hatsushika, Takanori Matsumoto, Nobuo Yamazaki, Hirotaka
  Kamahori, Kiyotoshi Takahashi, Shinji Kadokura, Koji Wada, Koji Kato, Ryo
  Oyama, Tomoaki Ose, Nobutaka Mannoji, and Ryusuke Taira}, {\em The jra-25
  reanalysis}, Journal of the Meteorological Society of Japan. Ser. II, 85
  (2007), pp.~369--432.

\bibitem{Packard_Crutchfield_Farmer_Shaw_1980}
{\sc N.~H. Packard, J.~P. Crutchfield, J.~D. Farmer, and R.~S. Shaw}, {\em
  Geometry from a time series}, Phys. Rev. Lett., 45 (1980), pp.~712--716.

\bibitem{Sauer_Yorke_Casdagli_1991}
{\sc T.~{Sauer}, J.~A. {Yorke}, and M.~{Casdagli}}, {\em {Embedology}}, Journal
  of Statistical Physics, 65 (1991), pp.~579--616.

\bibitem{Stephenson_Hannachi_ONeill_2004}
{\sc D.~B. {Stephenson}, A.~{Hannachi}, and A.~{O'Neill}}, {\em {On the
  existence of multiple climate regimes}}, Quarterly Journal of the Royal
  Meteorological Society, 130 (2004), pp.~583--605.

\bibitem{Takens_1981}
{\sc Floris Takens}, {\em Detecting strange attractors in turbulence}, in
  Dynamical systems and turbulence, vol.~898 of Lecture Notes in Mathematics,
  Springer, 1981, pp.~366--381.


\bibitem{Walker_1924}
{\sc G.~T. Walker}, {\em {Correlation in seasonal variation of weather. IX. 
A further study of world weather}}, Mem. Ind. Meteor. Dept., 24, (1924), pp.~275-333.

\bibitem{Wallace_Gutzler_1981}
{\sc J.~M.~Wallace and D.~S.~Gutzler}, {\em {Teleconnections in the geopotential height field 
during the Northern Hemisphere winter}}, Mon. Wea. Rev., 109 (1981), pp.~784-812.

\bibitem{Woolings_Pinto_Santos}
{\sc T. M. Woolings, J. G. Pinto, and J. A. Santos}, {\em Dynamical evolution of North Atlantic 
ridges and poleward jet stream displacements}. J. Atmos. Sci., 68 (2011), pp. 954--963.



\end{thebibliography}

\end{document}